\documentclass[12pt]{article}
\usepackage{amsmath,amsfonts,amssymb,amsthm}
\input amssym.def
\begin{document}
\def \Z{\Bbb Z}
\def \C{\Bbb C}
\def \R{\Bbb R}
\def \Q{\Bbb Q}
\def \N{\Bbb N}

\def \A{{\mathcal{A}}}
\def \D{{\mathcal{D}}}
\def \E{{\mathcal{E}}}
\def \E{{\mathcal{E}}}
\def \H{\mathcal{H}}
\def \S{{\mathcal{S}}}
\def \wt{{\rm wt}}
\def \tr{{\rm tr}}
\def \span{{\rm span}}
\def \Res{{\rm Res}}
\def \Der{{\rm Der}}
\def \End{{\rm End}}
\def \Ind {{\rm Ind}}
\def \Irr {{\rm Irr}}
\def \Aut{{\rm Aut}}
\def \GL{{\rm GL}}
\def \Hom{{\rm Hom}}
\def \mod{{\rm mod}}
\def \ann{{\rm Ann}}
\def \ad{{\rm ad}}
\def \rank{{\rm rank}\;}
\def \<{\langle}
\def \>{\rangle}

\def \g{{\frak{g}}}
\def \h{{\hbar}}
\def \k{{\frak{k}}}
\def \sl{{\frak{sl}}}
\def \gl{{\frak{gl}}}

\def \be{\begin{equation}\label}
\def \ee{\end{equation}}
\def \bex{\begin{example}\label}
\def \eex{\end{example}}
\def \bl{\begin{lem}\label}
\def \el{\end{lem}}
\def \bt{\begin{thm}\label}
\def \et{\end{thm}}
\def \bp{\begin{prop}\label}
\def \ep{\end{prop}}
\def \br{\begin{rem}\label}
\def \er{\end{rem}}
\def \bc{\begin{coro}\label}
\def \ec{\end{coro}}
\def \bd{\begin{de}\label}
\def \ed{\end{de}}

\newcommand{\m}{\bf m}
\newcommand{\n}{\bf n}
\newcommand{\nno}{\nonumber}
\newcommand{\nord}{\mbox{\scriptsize ${\circ\atop\circ}$}}
\newtheorem{thm}{Theorem}[section]
\newtheorem{prop}[thm]{Proposition}
\newtheorem{coro}[thm]{Corollary}
\newtheorem{conj}[thm]{Conjecture}
\newtheorem{example}[thm]{Example}
\newtheorem{lem}[thm]{Lemma}
\newtheorem{rem}[thm]{Remark}
\newtheorem{de}[thm]{Definition}
\newtheorem{hy}[thm]{Hypothesis}
\makeatletter \@addtoreset{equation}{section}
\def\theequation{\thesection.\arabic{equation}}
\makeatother \makeatletter

\begin{center}
{\Large \bf Twisted tensor products of nonlocal vertex algebras}
\end{center}

\begin{center}
{Haisheng Li \\
Department of Mathematical Sciences, Rutgers University, Camden, NJ 08102\\
Jiancai Sun\\
Department of Mathematics, Shanghai University, Shanghai 200444,
China}
\end{center}

\begin{abstract}
In this paper we introduce and study a twisted tensor product
construction of nonlocal vertex algebras. Among the main results, we
establish a universal property and give a characterization of a
twisted tensor product. Furthermore, we give a construction of
modules for a twisted tensor product. We also show that smash
products studied by one of us before can be realized as twisted
tensor products.
\end{abstract}

\section{Introduction}
For associative algebras, there is a notion of twisted tensor
product, generalizing the notion of tensor product. Let $A$ and $B$
be two associative unital algebras. A twisted tensor product
$A\otimes_{R}B$ (see \cite{csv}, \cite{vv}) is associated to each
twisting operator $R$ in the sense that $R: B\otimes A\rightarrow
A\otimes B$ is a linear map satisfying
\begin{eqnarray*}
&&R(1\otimes a)=a\otimes 1, \ \ \ R(b\otimes 1)=1\otimes b\ \ \
\mbox{ for }a\in A,\; b\in B,\\
&&R(\mu_{B}\otimes 1)=(1\otimes \mu_{B})R_{12}R_{23},\ \ \
R(1\otimes \mu_{A})=(\mu_{A}\otimes 1)R_{23}R_{12},
\end{eqnarray*}
where $\mu_{A}$ and $\mu_{B}$ denote the multiplications of $A$ and
$B$. The multiplication $\mu$ of $A\otimes_{R}B$ is given by
$$\mu (a\otimes b\otimes a'\otimes b')=(\mu_{A}\otimes \mu_{B})
(a\otimes R(b\otimes a')\otimes b')$$ for $a,a'\in A,\ b,b'\in B$.
 Twisted tensor product, as well as its generalizations,
provides an effective tool to construct new algebras for various
purposes (cf. \cite{ppo}).

Vertex algebras are both analogues and generalizations of
commutative and associative unital algebras, while nonlocal vertex
algebras (or field algebras in the sense of \cite{bk}) are analogues
and generalizations of associative unital algebras. In \cite{ek},
one of an important series of papers, Etingof and Kazhdan developed
a fundamental theory of quantum vertex operator algebras in the
sense of formal deformation, where quantum vertex operator algebras
are ($\hbar$-adic) nonlocal vertex algebras (over $\C[[\hbar]]$)
which satisfy what was called $\S$-locality. Partly motivated by
this, we developed a theory of (weak) quantum vertex algebras (see
\cite{li-qva1}, \cite{li-qva2}), where weak quantum vertex algebras
are generalizations of vertex superalgebras, instead of formal
deformations. Weak quantum vertex algebras in this sense are
nonlocal vertex algebras that satisfy a variation of
Etingof-Kazhdan's $\S$-locality. In this developing theory,
constructing interesting examples of quantum vertex algebras is one
of the most important problems.

In this paper, we study twisted tensor products of nonlocal vertex
algebras and of (weak) quantum vertex algebras. The main purpose is
to build various tools for constructing new interesting quantum
vertex algebras. Let $U$ and $V$ be two nonlocal vertex algebras. We
define a twisting operator to be a linear map
$$R(x):\ V\otimes U\rightarrow U\otimes V\otimes \C((x)),$$
satisfying a set of conditions which are stringy analogues of those
listed before for a twisting operator with associative algebras. The
underlying space of the twisted tensor product $U\otimes_{R}V$
associated to $R$ is $U\otimes V$, while the vacuum vector is  ${\bf
1}_{U}\otimes {\bf 1}_{V}$ and  the vertex operator map, denoted by
$Y_{R}$, is given by
$$Y_{R}(u\otimes v,x)(u'\otimes v')=(Y_{U}(x)\otimes Y_{V}(x))(u\otimes
R(x)(v\otimes u')\otimes v')$$ for $u,u'\in U,\ v,v'\in V$. It is
proved that $U\otimes_{R}V$ is a nonlocal vertex algebra, containing
$U$ and $V$ canonically as subalgebras which satisfy a certain
commutation relation. (If both $U$ and $V$ are weak quantum vertex
algebras, it is proved that $U\otimes_{R}V$ is a weak quantum vertex
algebra.)  On the other hand, it is proved that if a nonlocal vertex
algebra $K$, which is non-degenerate in the sense of \cite{ek},
contains subalgebras $U$ and $V$ satisfying a certain commutation
relation, then $K$ is isomorphic to the twisted tensor product
$U\otimes_{R}V$ with respect to a twisting operator $R(x)$. Also
established in this paper is a universal property for the twisted
tensor product $U\otimes_{R}V$, similar to the one for the ordinary
tensor product $U\otimes V$. Regarding $U\otimes_{R}V$-modules, it
is proved that a $U$-module structure and a $V$-module structure
compatible in a certain sense on a vector space $W$ give rise to a
$U\otimes_{R}V$-module structure canonically.

In \cite{li-smash}, a notion of nonlocal vertex bialgebra and a
notion of nonlocal vertex module-algebra for a nonlocal vertex
bialgebra were formulated, and a smash product construction of
nonlocal vertex algebras was established. Given a nonlocal vertex
bialgebra $H$ and for a nonlocal vertex $H$-module algebra $U$, we
have a smash product $U\sharp H$. In the present paper, we slightly
generalize the smash product construction with $H$ replaced by what
we call a nonlocal vertex $H$-comodule-algebra, and we show that the
smash product $U\sharp V$ is a twisted tensor product with respect
to a canonical twisting operator.

This paper was partly motivated by a recent study \cite{ls} on
regular representations for what we called M\"{o}bius quantum vertex
algebras.  Previously, a theory of regular representations for a
general vertex operator algebra $V$ was developed in
\cite{li-regular}, where the regular representation space was proved
to have a canonical module structure for the ordinary tensor product
$V\otimes V$. Furthermore, under suitable assumptions on $V$, a
result of Peter-Weyl type was obtained. In \cite{ls} we extended
this theory for a M\"{o}bius quantum vertex algebra $V$, and we
proved that the regular representation space has a canonical module
structure for a certain twisted tensor product of $V$ with $V$,
instead of the ordinary tensor product. This motivated us to study
more general twisted tensor products of nonlocal vertex algebras.

We mention that Anguelova and Bergvelt studied a notion of what they
called $H_{D}$-quantum vertex algebra in \cite{ab} and they gave a
construction by employing Borcherds' bicharacter construction (see
\cite{bor}).

This paper is organized as follows: In Section 2, we present some
basic notions and twisted tensor product. In Section 3, we show that
twisted tensor products are (weak) quantum vertex algebras if the
tensor factors are (weak) quantum vertex algebras. In Section 4, we
show that smash products of nonlocal vertex algebras over a nonlocal
vertex bialgebra are twisted tensor products.

{\bf Acknowledgment:} H. Li gratefully acknowledges financial
support from the United States National Security Agency under grant
H98230-11-1-0161. Part of this work was done during a visit by H. Li
at Kavli Institute for Theoretical Physics China, CAS, Beijing,
2010. We would like to thank KITP for the financial support during
the visit; This research was supported in part by the Project of
Knowledge Innovation Program (PKIP) of Chinese Academy of Sciences,
Grant No. KJCX2.YW.W10.

\section{Twisted tensor product nonlocal vertex algebras}
In this section, first we define the notion of twisting operator and
construct the twisted tensor product nonlocal vertex algebra. Then
we establish a universal property and give a characterization of the
twisted tensor product. We also present an existence theorem for a
module structure for the twisted tensor product.

We begin by recalling the notion of nonlocal vertex algebra. A {\em
nonlocal vertex algebra} (see \cite{li-g1}, cf. \cite{bk}) is a
vector space $V$, equipped with a linear map
\begin{eqnarray*}
Y(\cdot,x):&&V\rightarrow \Hom (V,V((x)))\subset (\End V)[[x,x^{-1}]]\\
&&v\mapsto Y(v,x)=\sum_{n\in\Z}v_{n}x^{-n-1}\ \ (\mbox{where
}v_{n}\in \End V),
\end{eqnarray*}
and equipped with a vector ${\bf 1}\in V$, satisfying the conditions
that for $v\in V,$
\begin{eqnarray}
Y({\bf 1},x)v=v,\ \ \ Y(v,x){\bf 1}\in V[[x]]\ \ \mbox{ and }\ \
\lim_{x\rightarrow 0}Y(v,x){\bf 1}=v,
\end{eqnarray}
and that for $u, v,w\in V$, there exists a nonnegative integer $k$
such that
\begin{eqnarray}\label{eweak-assoc}
(x_{0}+x_{2})^{k}Y(u,x_{0}+x_{2})Y(v,x_{2})w=(x_{0}+x_{2})^{k}Y(Y(u,x_{0})v,x_{2})w
\end{eqnarray}
({\em weak associativity}).

We sometimes denote a nonlocal vertex algebra by a triple $(V,Y,{\bf
1})$, to emphasize the {\em vertex operator map} $Y$ and the {\em
vacuum vector} ${\bf 1}$.

Let $V$ be a nonlocal vertex algebra. Define a linear operator $\D$
on $V$ by
\begin{eqnarray}
\D(v)=v_{-2}{\bf 1}\ \ \mbox{ for }v\in V.
\end{eqnarray}
Then
\begin{eqnarray}
[\D,Y(v,x)]=Y(\D v,x)=\frac{d}{dx}Y(v,x)\ \ \mbox{ for }v\in
V\label{d}.
\end{eqnarray}

For a nonlocal vertex algebra $V$, a {\em $V$-module} is a vector
space $W$, equipped with a linear map
$$Y_{W}(\cdot,x):\ V\rightarrow \Hom (W,W((x)))\subset (\End
W)[[x,x^{-1}]],$$ satisfying the conditions that $Y_{W}({\bf
1},x)=1_{W}$ (the identity operator on $W$) and that for $u,v\in V,\
w\in W$, there exists a nonnegative integer $l$ such that
\begin{eqnarray*}
(x_{0}+x_{2})^{l}Y_{W}(u,x_{0}+x_{2})Y_{W}(v,x_{2})w
=(x_{0}+x_{2})^{l}Y_{W}(Y(u,x_{0})v,x_{2})w.
\end{eqnarray*}

Let $(W,Y_{W})$ be a $V$-module. Note that for $u,v\in V$,
$$(Y_{W}(u,x_{1})Y_{W}(v,x_{2}))|_{x_{1}=x_{0}+x_{2}}=Y_{W}(u,x_{0}+x_{2})Y_{W}(v,x_{2}),$$
which by the formal Taylor theorem equals
$e^{x_{2}\frac{\partial}{\partial
x_{0}}}(Y_{W}(u,x_{0})Y_{W}(v,x_{2}))$, exists, but
$(Y_{W}(u,x_{1})Y_{W}(v,x_{2}))|_{x_{1}=x_{2}+x_{0}}$ in general
does {\em not} exist. On the other hand, if $A(x_{1},x_{2})$ is an
element of $\Hom (W,W((x_{1},x_{2})))$, then both
$A(x_{1},x_{2})|_{x_{1}=x_{2}+x_{0}}$ and
$A(x_{1},x_{2})|_{x_{1}=x_{0}+x_{2}}$ exist.

As we shall need, we recall another form of weak associativity from
\cite{ltw}.

\bl{ldef-module} Let $V$ be a nonlocal vertex algebra. In the
definition of a $V$-module, the weak associativity axiom can be
equivalently replaced with the property that for $u,v\in V$, there
exists a nonnegative integer $k$ such that
\begin{eqnarray}\label{emodule-compatibility}
(x_{1}-x_{2})^{k}Y_{W}(u,x_{1})Y_{W}(v,x_{2})\in \Hom
(W,W((x_{1},x_{2})))
\end{eqnarray}
and
\begin{eqnarray}\label{eweak-assoc-2}
\left((x_{1}-x_{2})^{k}Y_{W}(u,x_{1})Y_{W}(v,x_{2})\right)|_{x_{1}=x_{2}+x_{0}}
=x_{0}^{k}Y_{W}(Y(u,x_{0})v,x_{2}).
\end{eqnarray}
Furthermore, for a $V$-module $W$ and for any $u,v\in V$,
(\ref{eweak-assoc-2}) holds for every nonnegative integer $k$
satisfying (\ref{emodule-compatibility}). \el

\begin{proof} The first part is a special case of Lemma 2.9 of \cite{ltw} with $\sigma=1$.
As for the second part, suppose that $k$ is a nonnegative integer
such that both (\ref{emodule-compatibility}) and
(\ref{eweak-assoc-2}) hold, and let $k'$ be any nonnegative integer
such that (\ref{emodule-compatibility}) holds. Note that
$$\left((x_{1}-x_{2})^{r}Y_{W}(u,x_{1})Y_{W}(v,x_{2})\right)|_{x_{1}=x_{2}+x_{0}}$$
 for $r=k$ and for $r=k'$ both exist in $\Hom (W,((x_{2}))[[x_{0}]])$.
 Then
\begin{eqnarray*}
&&x_{0}^{k}\left((x_{1}-x_{2})^{k'}Y_{W}(u,x_{1})Y_{W}(v,x_{2})\right)|_{x_{1}=x_{2}+x_{0}}\\
&=&\left((x_{1}-x_{2})^{k}\right)|_{x_{1}=x_{2}+x_{0}}\cdot
\left((x_{1}-x_{2})^{k'}Y_{W}(u,x_{1})Y_{W}(v,x_{2})\right)|_{x_{1}=x_{2}+x_{0}}\\
&=&\left((x_{1}-x_{2})^{k+k'}Y_{W}(u,x_{1})Y_{W}(v,x_{2})\right)|_{x_{1}=x_{2}+x_{0}}\\
&=&\left((x_{1}-x_{2})^{k'}\right)|_{x_{1}=x_{2}+x_{0}}\cdot
\left((x_{1}-x_{2})^{k}Y_{W}(u,x_{1})Y_{W}(v,x_{2})\right)|_{x_{1}=x_{2}+x_{0}}\\
&=&x_{0}^{k+k'}Y_{W}(Y(u,x_{0})v,x_{2}),
\end{eqnarray*}
{}from which we immediately get
\begin{eqnarray*}
\left((x_{1}-x_{2})^{k'}Y_{W}(u,x_{1})Y_{W}(v,x_{2})\right)|_{x_{1}=x_{2}+x_{0}}
=x_{0}^{k'}Y_{W}(Y(u,x_{0})v,x_{2}),
\end{eqnarray*}
as desired.
\end{proof}

For a nonlocal vertex algebra $(V,Y,{\bf 1})$, we follow \cite{ek}
to define a linear map
$$Y(x):\ V\otimes V\rightarrow V((x))$$
by
$$Y(x)(u\otimes v)=Y(u,x)v\ \ \ \mbox{ for }u,v\in V.$$
Similarly, for a $V$-module $(W,Y_{W})$, we denote by $Y_{W}(x)$ the
associated linear map
$$Y_{W}(x): \ V\otimes W\rightarrow W((x)).$$

Let $U$ and $V$ be two nonlocal vertex algebras. We have an
(ordinary) tensor product nonlocal vertex algebra $U\otimes V$,
where the vacuum vector is ${\bf 1}\otimes {\bf 1}$ and the vertex
operator map is given by
$$Y(u\otimes v,x)(u'\otimes v')=Y(u,x)u'\otimes Y(v,x)v'
\ \ \ \mbox{ for }u,u'\in U,\ v,v'\in V.$$ That is,
$$Y_{U\otimes V}(x)=\left(Y_{U}(x)\otimes
Y_{V}(x)\right)\sigma^{23},$$ where $\sigma^{23}$ is the linear
operator on $(U\otimes V)^{\otimes 2}$, defined by
$$\sigma^{23}(u\otimes v\otimes u'\otimes v')=u\otimes u'\otimes
v\otimes v'$$
for $u,u'\in U,\ v,v'\in V.$

\bd{dtwisting} {\em Let $U$ and $V$ be nonlocal vertex algebras. A
{\em twisting operator} for the ordered pair $(U,V)$ is a linear map
$$R(x):\  V\otimes U\rightarrow U\otimes V\otimes \C((x)),$$
satisfying the following conditions:
\begin{eqnarray}
&&R(x)(v\otimes {\bf
1})={\bf 1}\otimes v\ \ \ \mbox{ for }v\in V,\\
&&R(x)({\bf 1}\otimes u)=u\otimes {\bf 1}\ \ \ \mbox{ for }u\in U,\\
&&R(x_{1})(1\otimes Y(x_{2}))=(Y(x_{2})\otimes
1)R^{23}(x_{1})R^{12}(x_{1}+x_{2}),\label{sy12}\\
&&R(x_{1})(Y(x_{2})\otimes 1)=(1\otimes
Y(x_{2}))R^{12}(x_{1}-x_{2})R^{23}(x_{1}).\label{ry12}
\end{eqnarray}} \ed

We say that a twisting operator $R(x)$ is {\em invertible} if
$R(x)$, viewed as a $\C((x))$-linear map from $V\otimes U\otimes
\C((x))$ to $U\otimes V\otimes \C((x))$, is invertible. The inverse
of an invertible $R(x)$ is a $\C((x))$-linear map $R^{-1}(x)$ {}from
$U\otimes V\otimes \C((x))$ to $V\otimes U\otimes \C((x))$. We often
consider $R^{-1}(x)$ as a $\C$-linear map
$$R^{-1}(x):\ \ U\otimes V\rightarrow V\otimes U\otimes \C((x)).$$

The following is straightforward:

\bl{linverse} If $R(x)$ is an invertible twisting operator for the
ordered pair $(U,V)$, then  $R^{-1}(-x)$ is an invertible twisting
operator for the ordered pair $(V,U)$. \el

We now present the twisted tensor product.

\bt{twqva} Let $U,V$ be nonlocal vertex algebras and let $R(x)$ be a
twisting operator of the ordered pair $(U,V)$. Set
\begin{eqnarray}
Y_{R}(x)=(Y(x)\otimes Y(x))R^{23}(-x).\label{pgloble}
\end{eqnarray}
Then $(U\otimes V, Y_{R}, {\bf 1}\otimes {\bf 1})$ carries the
structure of a nonlocal vertex algebra, which contains $U$ and $V$
canonically as nonlocal vertex subalgebras. \et

\begin{proof} For $u, u'\in U,\ v, v' \in V$, by definition we have
\begin{eqnarray}
&&Y_{R}(x)(u\otimes v\otimes u'\otimes
v')=\sum_{i=1}^{r}f_{i}(-x)Y(u,x)u'^{(i)}\otimes
Y(v^{(i)},x)v',\label{pvtzv1}
\end{eqnarray}
where
\begin{eqnarray*}
R(x)(v\otimes u')=\sum_{i=1}^{r}u'^{(i)}\otimes v^{(i)}\otimes
f_{i}(x)\in U\otimes V\otimes \C((x)).
\end{eqnarray*}
As $$f_{i}(-x)\in \C((x)),\ \ Y(u,x)u'^{(i)}\in U((x)),\ \
Y(v^{(i)},x)v' \in V((x))$$ for $1\le i\le r$, we see that
$Y_{R}(u\otimes v,x)(u'\otimes v')$ exists in $(U\otimes V)((x))$.

For $u\in U,\ v\in V$, with $R(x)({\bf 1}\otimes u)=u\otimes {\bf
1}$, we have
\begin{eqnarray*}
Y_{R}({\bf 1}\otimes {\bf 1}, x)(u\otimes v) &=&(Y(x)\otimes
Y(x))({\bf 1}\otimes u\otimes {\bf
1}\otimes v)\\
&=&Y({\bf 1},x)u\otimes Y({\bf 1},x)v\\
&=&u\otimes v.
\end{eqnarray*}
On the other hand, for $u\in U,\; v\in V$, with $R(x)(v\otimes {\bf
1})={\bf 1}\otimes v$, we have
\begin{eqnarray}
Y_{R}(u\otimes v,x)({\bf 1}\otimes {\bf 1}) &=&(Y(x)\otimes
Y(x))R^{23}(-x)(u\otimes v\otimes {\bf
1}\otimes {\bf 1})\nonumber\\
&=&(Y(x)\otimes Y(x))(u\otimes {\bf
1}\otimes v\otimes {\bf 1})\nonumber\\
&=&Y(u,x){\bf 1}\otimes Y(v,x){\bf 1},
\end{eqnarray}
which lies in $(U\otimes V)[[x]]$, and
\begin{eqnarray*}
&&\lim_{x\rightarrow 0}Y_{R}(u\otimes v,x)({\bf 1}\otimes {\bf
1})=\lim_{x\rightarrow 0}Y(u,x){\bf 1}\otimes Y(v,x){\bf 1}=u\otimes
v.
\end{eqnarray*}

To see weak associativity, let $u, u', u''\in U,\ v, v', v'' \in V$.
Using (\ref{sy12}) we get
\begin{eqnarray}
&&Y_{R}(u\otimes v,x_{0}+x_{2})Y_{R}(u'\otimes v',x_{2})(u''\otimes v'')\nonumber\\
&=&(Y(x_{0}+x_{2})\otimes
Y(x_{0}+x_{2}))R^{23}(-x_{0}-x_{2})(1\otimes 1\otimes Y(x_{2})\otimes Y(x_{2}))\nonumber\\
&&\cdot R^{45}(-x_{2})(u\otimes v\otimes u'\otimes
v'\otimes u''\otimes v'')\nonumber\\
&=&(Y(x_{0}+x_{2})\otimes
Y(x_{0}+x_{2}))(1\otimes Y(x_{2})\otimes 1\otimes Y(x_{2}))R^{34}(-x_{0}-x_{2})\nonumber\\
&&\cdot R^{23}(-x_{0})R^{45}(-x_{2})(u\otimes v\otimes u'\otimes
v'\otimes u''\otimes v'').
\end{eqnarray}
On the other hand, using (\ref{ry12}) we get
\begin{eqnarray}
&&Y_{R}(Y_{R}(u\otimes v,x_{0})(u'\otimes v'),x_{2})(u''\otimes v'')\nonumber\\
&=&(Y(x_{2})\otimes
Y(x_{2}))R^{23}(-x_{2})(Y(x_{0})\otimes Y(x_{0})\otimes 1\otimes 1)\nonumber\\
&&\cdot R^{23}(-x_{0})(u\otimes v\otimes u'\otimes
v'\otimes u''\otimes v'')\nonumber\\
&=&(Y(x_{2})\otimes
Y(x_{2}))(Y(x_{0})\otimes 1\otimes Y(x_{0})\otimes 1)R^{34}(-x_{2}-x_{0})R^{45}(-x_{2})\nonumber\\
&&\cdot R^{23}(-x_{0})(u\otimes
v\otimes u'\otimes v'\otimes u''\otimes v'')\nonumber\\
&=&(Y(x_{2})\otimes Y(x_{2}))(Y(x_{0})\otimes 1\otimes
Y(x_{0})\otimes
1)R^{34}(-x_{2}-x_{0})R^{23}(-x_{0})\nonumber\\
&&\cdot R^{45}(-x_{2})(u\otimes v\otimes u'\otimes v'\otimes
u''\otimes v'').
\end{eqnarray}
Then the desired weak associativity relation follows. Thus
$(U\otimes V,Y_{R},{\bf 1}\otimes {\bf 1})$ carries the structure of
a nonlocal vertex algebra.

Furthermore, for $u,u'\in U$, as $R(x)({\bf 1}\otimes u')=u'\otimes
{\bf 1}$, we have
$$Y_{R}(u\otimes {\bf 1},x)(u'\otimes {\bf 1})=Y(u,x)u'\otimes
Y({\bf 1},x){\bf 1}=Y(u,x)u'\otimes {\bf 1}.$$ It follows that the
map $u\in U\mapsto u\otimes {\bf 1}\in U\otimes  V$ is a one-to-one
homomorphism of nonlocal vertex algebras. Similarly, the map $v\in
V\mapsto {\bf 1}\otimes v\in U\otimes V$ is a one-to-one
homomorphism of nonlocal vertex algebras. This concludes the proof.
\end{proof}

Denote by $U\otimes _{R}V$ the nonlocal vertex algebra obtained in
Theorem \ref{twqva}. We identify each element $u$ of $U$ with the
element $u\otimes {\bf 1}$ of $U\otimes_{R}V$ and identify each
element $v$ of $V$ with ${\bf 1}\otimes v$ of $U\otimes_{R}V$. For
$u,u'\in U,\; v,v'\in V$, we have
\begin{eqnarray}
&&Y_{R}(u,x)(u'\otimes v') =Y(u,x)u'\otimes v',\label{eu-action}\\
&&Y_{R}(v,x)(u'\otimes v')=(1\otimes Y(x))R^{12}(-x)(v\otimes
u'\otimes v').
\end{eqnarray}

The following two propositions give more information about
$U\otimes_{R}V$:

\bp{pmoreproperty}  The $\D$-operator of $U\otimes_{R}V$ is given by
\begin{eqnarray}
\D_{U\otimes_{R}V}=\D\otimes 1+1\otimes \D,
\end{eqnarray}
where the two $\D$'s denote the $\D$-operators of $U$ and $V$,
respectively. Furthermore, we have
\begin{eqnarray}
&&Y_{R}(u,x)v\in (U\otimes V)[[x]],\label{euv-regularity}\\
&&Y_{R}(v,x)u=e^{x(\D\otimes 1+1\otimes\D)}Y_{R}(-x)R(-x)(v\otimes
u)\label{evu-symmetry}
\end{eqnarray}
for $u\in U,\ v\in V$. \ep

\begin{proof} Let $u\in U,\ v\in V$. {}From definition we have
\begin{eqnarray*}
&&\Res_{x}x^{-2}Y_{R}(u\otimes v,x)({\bf 1}\otimes {\bf 1})\nonumber\\
&=&\Res_{x}x^{-2}(Y(u,x){\bf 1}\otimes Y(v,x){\bf 1})\nonumber\\
&=&\Res_{x}x^{-2}e^{x(\D\otimes 1+1\otimes \D)}(u\otimes
v)\nonumber\\
&=&(\D\otimes 1+1\otimes\D)(u\otimes v).\label{dvv}
\end{eqnarray*}
Thus the $\D$-operator is given by $\D\otimes 1+1\otimes \D$.

Let $u\in U,\ v\in V$. As
$$R^{23}(-x)(u\otimes {\bf 1}\otimes {\bf 1}\otimes v)
=u\otimes {\bf 1}\otimes {\bf 1}\otimes v,$$ we have
\begin{eqnarray*}
Y_{R}(u,x)v=Y_{R}(u\otimes {\bf 1},x)({\bf 1}\otimes v)=Y(u,x){\bf
1}\otimes Y({\bf 1},x)v=Y(u,x){\bf 1}\otimes v,
\end{eqnarray*}
which implies (\ref{euv-regularity}) and
\begin{eqnarray}
u\otimes v=u_{-1}{\bf 1}\otimes v=(u\otimes {\bf 1})_{-1}({\bf
1}\otimes v).
\end{eqnarray}

For $u\in U,\ v\in V$, we have
\begin{eqnarray*}
&&Y_{R}({\bf 1}\otimes v,x)(u\otimes {\bf 1})\\
&=&(Y(x)\otimes Y(x))R^{23}(-x)({\bf 1}\otimes v\otimes u\otimes {\bf 1})\\
&=&e^{x(\D\otimes 1+1\otimes \D)}(Y(-x)\otimes
Y(-x))R^{14}(-x)(v\otimes {\bf 1}\otimes {\bf 1}\otimes u)\\
&=&e^{x(\D\otimes 1+1\otimes \D)}Y_{R}(-x)R^{14}(-x)(v\otimes {\bf
1}\otimes {\bf 1}\otimes u),
\end{eqnarray*}
as $R(-x)({\bf 1}\otimes {\bf 1})={\bf 1}\otimes {\bf 1}$. This
proves (\ref{evu-symmetry}).
\end{proof}

\bp{prelations-prod} Let $U,V, R(x)$ be given as in Theorem
\ref{twqva}. For $u\in U,\ v\in V$, there exists a nonnegative
integer $k$ such that
\begin{eqnarray}
&&(x_{1}-x_{2})^{k}Y_{R}(v,x_{1})Y_{R}(u,x_{2})w\nonumber\\
&=&(x_{1}-x_{2})^{k} Y_{R}(x_{2})(1\otimes
Y_{R}(x_{1}))R^{12}(x_{2}-x_{1})(v\otimes u\otimes w)
\end{eqnarray}
for every $w\in U\otimes_{R}V$. Furthermore, if $R(x)$ is
invertible, we also have
\begin{eqnarray}
&&Y_{R}(u,x)v=e^{x(\D\otimes 1+1\otimes \D)}Y_{R}(-x)R^{-1}(x)(u\otimes v),\label{eskew-symmetry-1}\\
&&Y_{R}(u,x_{1})Y_{R}(v,x_{2})w=Y_{R}(x_{2})(1\otimes
Y_{R}(x_{1}))(R^{-1})^{12}(-x_{2}+x_{1})(u\otimes v\otimes w)\ \ \ \
\ \label{eyruv-comm}
\end{eqnarray}
for $u\in U,\ v\in V$ and for $w\in U\otimes_{R}V$. \ep

\begin{proof} For $u\in U,\ v\in V$, skew-symmetry relation (\ref{evu-symmetry}) holds.
{}From \cite{li-qva1} (Proposition 5.2), there exists a nonnegative
integer $k$ such that
\begin{eqnarray}
&&(x_{1}-x_{2})^{k}Y_{R}({\bf 1}\otimes v,x_{1})Y_{R}(u\otimes {\bf
1},x_{2})w\nonumber\\
&=&(x_{1}-x_{2})^{k}Y_{R}(x_{2})(1\otimes
Y_{R}(x_{1}))R^{14}(x_{2}-x_{1})(v\otimes {\bf 1}\otimes {\bf
1}\otimes u\otimes w)
\end{eqnarray}
for every $w\in U\otimes_{R}V$.

Now we assume that $R(x)$ is invertible. For $u\in U,\ v\in V$, we
have
\begin{eqnarray*}
&&Y_{R}(u,x)v\ \left(=Y_{R}(u\otimes {\bf 1},x)({\bf 1}\otimes v)\right)\\
&=&(Y(x)\otimes
Y(x))R^{23}(-x)(u\otimes {\bf 1}\otimes {\bf 1}\otimes v)\\
&=&Y(u,x){\bf 1}\otimes Y({\bf 1},x)v\\
&=&e^{x\D}Y({\bf 1},-x)u\otimes e^{x\D}Y(v,-x){\bf 1}\\
&=&e^{x(\D\otimes 1+1\otimes \D)}(Y(-x)\otimes Y(-x))({\bf 1}\otimes
u\otimes v\otimes {\bf 1})\\
&=&e^{x(\D\otimes 1+1\otimes \D)}Y_{R}(-x)(R^{-1})^{23}(x)({\bf
1}\otimes u\otimes v\otimes {\bf 1})\\
&=&e^{x(\D\otimes 1+1\otimes \D)}Y_{R}(-x)R^{-1}(x)(u\otimes v),
\end{eqnarray*}
proving (\ref{eskew-symmetry-1}). Again, from \cite{li-qva1}
(Proposition 5.2), there exists a nonnegative integer $k$ such that
\begin{eqnarray*}
&&(x_{1}-x_{2})^{k}Y_{R}(u,x_{1})Y_{R}(v,x_{2})w\\
 &=&(x_{1}-x_{2})^{k}Y_{R}(x_{2})(1\otimes
Y_{R}(x_{1}))(R^{-1})^{12}(-x_{2}+x_{1})(u\otimes v\otimes w)
\end{eqnarray*}
for every $w\in U\otimes_{R}V$. Combining this with weak
associativity we get
\begin{eqnarray*}
&&x_{0}^{-1}\delta\left(\frac{x_{1}-x_{2}}{x_{0}}\right)Y_{R}(u,x_{1})Y_{R}(v,x_{2})w
\nonumber\\
&&\hspace{0.5cm}-x_{0}^{-1}\delta\left(\frac{x_{2}-x_{1}}{-x_{0}}\right)Y_{R}(x_{2})(1\otimes
Y_{R}(x_{1}))(R^{-1})^{12}(-x_{2}+x_{1})(u\otimes v\otimes w)\nonumber\\
&=&x_{2}^{-1}\delta\left(\frac{x_{1}-x_{0}}{x_{2}}\right)Y_{R}(Y_{R}(u,x_{0})v,x_{2})w.
\end{eqnarray*}
Note that $Y_{R}(u,x_{0})v\in (U\otimes V)[[x_{0}]]$ by
(\ref{euv-regularity}). Applying $\Res_{x_{0}}$ to the above Jacobi
identity we obtain (\ref{eyruv-comm}).
\end{proof}

The twisted tensor product $U\otimes_{R}V$ has a universal property
just as the ordinary tensor product $U\otimes V$ does (cf.
\cite{fhl}, \cite{ll}).

\bt{tuniversal} Let $U, V$ be nonlocal vertex algebras and let
$R(x)$ be a twisting operator for $(U,V)$. Let $K$ be any nonlocal
vertex algebra and let $\psi_{1}: U\rightarrow K$,  $\psi_{2}:
V\rightarrow K$ be any homomorphisms, satisfying the condition that
for $u\in U,\ v\in V$,
\begin{eqnarray}
&&Y(\psi_{1}(u),x)\psi_{2}(v)\in K[[x]],\\
&&Y(\psi_{2}(v),x)\psi_{1}(u)=e^{x\D}Y(-x)(\psi_{1}\otimes
\psi_{2})R(-x)(v\otimes u).\label{euniversal-skew}
\end{eqnarray}
Then the linear
map $\psi: U\otimes_{R}V\rightarrow K$, defined by
\begin{eqnarray*}
\psi(u\otimes v)= \psi_{1}(u)_{-1}\psi_{2}(v) \ \ \ \mbox{ for }u\in
U,\ v\in V,
\end{eqnarray*}
is a homomorphism of nonlocal vertex algebras, which extends
$\psi_{1}$ and $\psi_{2}$ uniquely.
 \et

\begin{proof} For $u\in U,\ v\in V$, as $Y(\psi_{1}(u),x)\psi_{2}(v)\in
K[[x]]$ by assumption, we have
$$\psi(u\otimes v)=\Res_{x}x^{-1}Y(\psi_{1}(u),x)\psi_{2}(v)
=\lim_{x\rightarrow 0}Y(\psi_{1}(u),x)\psi_{2}(v).$$ It is clear
that linear map $\psi$ extends both $\psi_{1}$ and $\psi_{2}$. It is
also clear that $\psi({\bf 1}\otimes {\bf 1})={\bf 1}$. To prove
that $\psi$ is a homomorphism of nonlocal vertex algebras, we must
prove
\begin{eqnarray}\label{ehomomorphism}
\psi(Y_{R}(u\otimes v,x)(u'\otimes v'))=Y(\psi(u\otimes
v),x)\psi(u'\otimes v')\ \ \mbox{ for }u,u'\in U,\ v,v'\in V.
\end{eqnarray}
Through the homomorphisms $\psi_{1}$ and $\psi_{2}$, $K$ becomes a
$U$-module and a $V$-module. Let $v\in V$. {}From assumption, we
have $Y(\psi_{1}(u),x)\psi_{2}(v)\in K[[x]]$ for every $u\in U.$ By
a result of \cite{li-qva1} (Lemma 6.1), for every fixed $v\in V$,
the map $u\in U\mapsto \psi_{1}(u)_{-1}\psi_{2}(v)\in K$ is a
$U$-module homomorphism. Then, for $u',u\in U,\ v\in V$, we have
\begin{eqnarray*}
&&\psi( Y_{R}(u',x_{0})(u\otimes v))=\psi(Y(u',x_{0})u\otimes v)
=\lim_{x\rightarrow
0}Y(\psi_{1}Y(u',x_{0})u,x)\psi_{2}(v)\\
&=&Y(u',x_{0})\left(\psi_{1}(u)_{-1}\psi_{2}(v)\right)
=Y(\psi_{1}(u'),x_{0})\psi(u\otimes v)=Y(\psi(u'\otimes {\bf
1}),x_{0})\psi(u\otimes v).
\end{eqnarray*}
 This shows that
(\ref{ehomomorphism}) holds with $v={\bf 1}$. We next show that
(\ref{ehomomorphism}) also holds with $u={\bf 1}$. We have
\begin{eqnarray}\label{efirst}
&&\psi(Y_{R}({\bf 1}\otimes v,x)(u'\otimes v'))\nonumber \\
&=&\psi\left((1\otimes Y(x))R^{12}(-x)(v\otimes u'\otimes v')\right)\nonumber\\
&=&\lim_{x_{2}\rightarrow 0}Y(x_{2})(1\otimes Y(x))
(\psi_{1}\otimes\psi_{2}\otimes\psi_{2}) R^{12}(-x)(v\otimes
u'\otimes v').
\end{eqnarray}
With the assumption (\ref{euniversal-skew}), from \cite{li-qva1}
(Proposition 5.2), there exists a nonnegative integer $k$ such that
\begin{eqnarray*}
&&(x_{2}-x)^{k}Y(\psi_{2}(v),x)Y(\psi_{1}(u'),x_{2})\psi_{2}(v')\\
&=&(x_{2}-x)^{k}Y(x_{2})(1\otimes Y(x))
(\psi_{1}\otimes\psi_{2}\otimes\psi_{2})R^{12}(x_{2}-x)(v\otimes
u'\otimes v').
\end{eqnarray*}
Noticing that $R(x)(v\otimes u')\in U\otimes V\otimes \C((x))$, we
may replace $k$ with a bigger integer so that $x^{k}R(x)(v\otimes
u')\in U\otimes V\otimes \C[[x]]$ also holds. Then
\begin{eqnarray*}
&&(-x)^{k}Y(\psi({\bf 1}\otimes v),x)\psi(u'\otimes v')\\
&=&\lim_{x_{2}\rightarrow 0}(x_{2}-x)^{k}Y(\psi_{2}(v),x)Y(\psi_{1}(u'),x_{2})\psi_{2}(v')\\
&=&\lim_{x_{2}\rightarrow 0}(x_{2}-x)^{k}Y(x_{2})(1\otimes Y(x))
(\psi_{1}\otimes\psi_{2}\otimes\psi_{2})R^{12}(x_{2}-x)(v\otimes
u'\otimes v')\\
&=&(-x)^{k}\lim_{x_{2}\rightarrow 0}Y(x_{2})(1\otimes Y(x))
(\psi_{1}\otimes\psi_{2}\otimes\psi_{2})R^{12}(-x)(v\otimes
u'\otimes v').
\end{eqnarray*}
Thus
\begin{eqnarray*}
&&Y(\psi({\bf 1}\otimes v),x)\psi(u'\otimes v')\\
&=&\lim_{x_{2}\rightarrow 0}Y(x_{2})(1\otimes Y(x))
(\psi_{1}\otimes\psi_{2}\otimes\psi_{2})R^{12}(-x)(v\otimes
u'\otimes v').
\end{eqnarray*}
Combining this with (\ref{efirst}) we obtain
$$\psi(Y_{R}({\bf 1}\otimes v,x)(u'\otimes
v'))=Y(\psi({\bf 1}\otimes v),x)\psi(u'\otimes v'),$$ proving that
(\ref{ehomomorphism}) holds with $u={\bf 1}$. Since $U\otimes_{R}V$
as a nonlocal vertex algebra is generated by the subset $U\cup V$,
it follows that $\psi$ is a homomorphism of nonlocal vertex
algebras. The uniqueness assertion is clear as $u\otimes v=(u\otimes
{\bf 1})_{-1}({\bf 1}\otimes v)$ for $u\in U,\; v\in V$.
\end{proof}

The following is a characterization of $U\otimes_{R}V$ in terms of
$U$, $V$ and $R(x)$:

\bp{pcharacterization} Let $U,V$ and $R(x)$ be given as in Theorem
\ref{tuniversal}, and let $K$ be a nonlocal vertex algebra which
contains $U$ and $V$ as subalgebras, satisfying
\begin{eqnarray*}
&&Y(u,x)v\in K[[x]],\\
&&Y(v,x)u=e^{x\D}Y(-x)R(-x)(v\otimes u)\ \ \ \mbox{ for }u\in U,\
v\in V.
\end{eqnarray*}
Assume that $K$ as a nonlocal vertex algebra is generated by $U\cup
V$ and that $U$ as a $U$-module is irreducible and of countable
dimension (over $\C$). Then the linear map $\theta:
U\otimes_{R}V\rightarrow K$, defined by $\theta(u\otimes v)=u_{-1}v$
for $u\in U,\ v\in V$, is a nonlocal-vertex-algebras isomorphism.
\ep

\begin{proof} It follows from Theorem
\ref{tuniversal} that $\theta$ is a nonlocal-vertex-algebra
homomorphism. Now we prove that $\theta$ is a bijection. Since
$U\otimes_{R}V$ as a nonlocal vertex algebra is generated by $U\cup
V$ and since $K$ is also generated by $U\cup V$, it follows that
$\theta$ is onto. On the other hand, as $\theta|_{U}=1$, $\theta$ is
a $U$-module homomorphism. Consequently, $\ker \theta$ is a
$U$-submodule of $U\otimes_{R}V$. From (\ref{eu-action}), for any
subspace $P$ of $V$, $U\otimes P$ is a $U$-submodule. As $U$ as a
$U$-module is irreducible and of countable dimension (over $\C$), by
a version of Schur lemma we have $\End_{U}U=\C$.  It follows that
$\ker \theta=U\otimes B$ for some subspace $B$ of $V$. For any $b\in
B$, we have ${\bf 1}\otimes b\in \ker \theta$, so that
$b=\theta({\bf 1}\otimes b)=0$. Thus $B=0$ and $\ker \theta=0$,
proving that $\theta$ is injective.
\end{proof}

Recall from Lemma \ref{linverse} that for any invertible twisting
operator $R(x)$ for $(U,V)$, $R^{-1}(-x)$ is an invertible twisting
operator for $(V,U)$. Furthermore, we have:

\bp{pinverse} Let $R(x)$ be an invertible twisting operator for
$(U,V)$ such that
\begin{eqnarray}\label{assumption-inverse}
R(x)(v\otimes u)\in U\otimes V\otimes \C[[x]],\ \ \
R^{-1}(x)(u\otimes v)\in V\otimes U\otimes \C[[x]]
\end{eqnarray}
for $u\in U,\ v\in V$. Then the linear map $\psi:\
V\otimes_{R^{-1}(-x)}U\rightarrow U\otimes_{R}V$, defined by
$$\psi(v\otimes u)=v_{-1}u\ \ (\mbox{in }U\otimes_{R}V)\ \ \ \mbox{ for }v\in V,\ u\in U,$$
is a nonlocal-vertex-algebra isomorphism. \ep

\begin{proof} With the assumption (\ref{assumption-inverse}),
combining (\ref{evu-symmetry}) with (\ref{euv-regularity}) we get
\begin{eqnarray}\label{eyrvu}
Y_{R}(v,x)u\in (U\otimes_{R}V)[[x]]\ \ \mbox{ for }u\in U,\ v\in V.
\end{eqnarray}
{}From Proposition \ref{prelations-prod} we also have
$$Y_{R}(u,x)v=e^{x(\D\otimes 1+1\otimes
\D)}Y_{R}(-x)R^{-1}(x)(u\otimes v).$$ By Theorem \ref{tuniversal},
$\psi$ is a nonlocal-vertex-algebra homomorphism. Clearly,
$\psi|_{U}=1$ and $\psi_{V}=1$. On the other hand, consider
$V\otimes_{R^{-1}(-x)}U$ and denote $Y_{R^{-1}(-x)}$ simply by
$Y_{R^{-1}}$. For $u\in U,\ v\in V$, by (\ref{euv-regularity}) and
(\ref{evu-symmetry}) we have
\begin{eqnarray*}
&&Y_{R^{-1}}(v,x)u\in (V\otimes U)[[x]],\\
&&Y_{R^{-1}}(u,x)v=e^{x(\D\otimes 1+1\otimes
\D)}Y_{R^{-1}}(-x)R^{-1}(x)(u\otimes v).
\end{eqnarray*}
Combining these with (\ref{assumption-inverse}) we get
\begin{eqnarray*}
Y_{R^{-1}}(u,x)v\in (V\otimes U)[[x]],
\end{eqnarray*}
while from Proposition \ref{prelations-prod} we have
\begin{eqnarray*}
Y_{R^{-1}}(v,x)u=e^{x(\D\otimes 1+1\otimes
\D)}Y_{R^{-1}}(-x)R(-x)(v\otimes u).
\end{eqnarray*}
By Theorem \ref{tuniversal}, there is a nonlocal-vertex-algebra
homomorphism $\phi: U\otimes_{R}V\rightarrow V\otimes_{R^{-1}(-x)}U$
such that $\phi(u\otimes v)=u_{-1}v$ for $u\in U,\ v\in V$. Because
$\psi\circ \phi$ and $\phi\circ \psi$ are nonlocal-vertex-algebra
homomorphisms preserving both $U$ and $V$ element-wise, it follows
that $\psi\circ\phi=1$ and $\phi\circ \psi=1$. Therefore, $\psi$ is
a nonlocal-vertex-algebra isomorphism.
\end{proof}

Next, we shall establish a refinement of Proposition
\ref{pcharacterization}. As we need, we recall an important notion
due to Etingof and Kazhdan. A nonlocal vertex algebra $V$ is said to
be {\em non-degenerate} (see \cite{ek}) if for every positive
integer $n$, the linear map
$$Z_{n}: V^{\otimes n}\otimes \C((x_{1}))\cdots ((x_{n}))\rightarrow
V((x_{1}))\cdots ((x_{n})),$$ defined by
$$Z_{n}(v^{(1)}\otimes \cdots\otimes v^{(n)}\otimes f)
=fY(v^{(1)},x_{1})\cdots Y(v^{(n)},x_{n}){\bf 1}$$ for
$v^{(1)},\dots, v^{(n)}\in V,\; f\in \C((x_{1}))\cdots ((x_{n}))$,
is injective.

Following \cite{li-qva1}, for every positive integer $n$, we define
a linear map
$$\pi_{n}: V^{\otimes n}\otimes \C((x_{1}))\cdots ((x_{n}))\rightarrow
\Hom (V,V((x_{1}))\cdots ((x_{n})))$$ by
$$\pi_{n}(v^{(1)}\otimes \cdots\otimes v^{(n)}\otimes f)(w)
=fY(v^{(1)},x_{1})\cdots Y(v^{(n)},x_{n})w$$ for $v^{(1)},\dots,
v^{(n)}\in V,\; f\in \C((x_{1}))\cdots ((x_{n}))$ and for $w\in V$.
By definition we have
$$Z_{n}(v^{(1)}\otimes \cdots\otimes v^{(n)}\otimes f)
=\pi_{n}(v^{(1)}\otimes \cdots\otimes v^{(n)}\otimes f)({\bf 1}).
$$ We see that $\pi_{n}$ is injective if $Z_{n}$ is injective. In particular,
non-degeneracy implies that all the linear maps $\pi_{n}$ for $n\ge
1$ are injective.

For convenience we recall a notion from \cite{li-qva1}. Let $P$ be a
vector space and let $r$ be a positive integer.  For $A,B\in
P[[x_{1}^{\pm 1},x_{2}^{\pm 1},\dots,x_{r}^{\pm 1}]]$, we define
$A\sim B$ if there exists a nonzero polynomial
$p(x_{1},\dots,x_{r})$ such that
 $$p(x_{1},\dots,x_{r})A(x_{1},\dots,x_{r})=p(x_{1},\dots,x_{r})B(x_{1},\dots,x_{r}).$$
This is an equivalence relation. Furthermore, when restricted onto
the subspace $P((x_{r}))\cdots ((x_{2}))((x_{1}))$, this equivalence
relation becomes equality relation.

Just as with bialgebras, a braided tensor product nonlocal vertex
algebra arises whenever a nonlocal vertex algebra contains two
compatible subalgebras. The following can be considered as an
analogue of a result in \cite{majid}:

\bt{tcharacter} Let $K$ be a nonlocal vertex algebra that contains
subalgebras $U$ and $V$, satisfying the condition that for $u\in U,\
v\in V$,
\begin{eqnarray}
Y(u,x)v\in K[[x]]
\end{eqnarray}
and  there exist
$$u^{(i)}\in U,\ v^{(i)}\in V,\ f_{i}(x)\in \C((x))\ \
(i=1,\dots,r)$$ such that
\begin{eqnarray}\label{evu-commutator}
(x_{1}-x_{2})^{k}Y(v,x_{1})Y(u,x_{2})=(x_{1}-x_{2})^{k}
\sum_{i=1}^{r}f_{i}(x_{2}-x_{1})Y(u^{(i)},x_{2})Y(v^{(i)},x_{1})
\end{eqnarray}
for some nonnegative integer $k$. Assume that $K$ is non-degenerate
and that $K$ is generated by $U\cup V$. Then there exists a linear
map $R(x): V\otimes U\rightarrow U\otimes V\otimes \C((x))$, which
is uniquely determined by the condition that
\begin{eqnarray}\label{evuw-condition}
Y(v,x_{1})Y(u,x_{2})w\sim Y(x_{2})(1\otimes
Y(x_{1}))R^{12}(x_{2}-x_{1})(v\otimes u\otimes w)
\end{eqnarray}
for $u\in U,\ v\in V,\ w\in K$, and $R(x)$ is a twisting operator.
Furthermore, the linear map $\psi: U\otimes V\rightarrow K$, defined
by
$$\psi(u\otimes v)= u_{-1}v\ \ \ \mbox{ for }u\in U,\ v\in V,$$
is a nonlocal-vertex-algebra isomorphism from $U\otimes_{R}V$ onto
$K$.\et

\begin{proof} From the assumption, there exists a linear map
$R(x):V\otimes U\rightarrow U\otimes V\otimes \C((x))$, satisfying
(\ref{evuw-condition}). Suppose that $T(x)$ is another such linear
map. For $u\in U,\ v\in V,\ w\in K$, we have
\begin{eqnarray*}
Y(x_{2})(1\otimes Y(x_{1}))R^{12}(x_{2}-x_{1})(v\otimes u\otimes
w)\sim Y(x_{2})(1\otimes Y(x_{1}))T^{12}(x_{2}-x_{1})(v\otimes
u\otimes w).
\end{eqnarray*}
As the expressions on both sides lie in $(U\otimes
V)((x_{2}))((x_{1}))$, we must have
\begin{eqnarray*}
Y(x_{2})(1\otimes Y(x_{1}))R^{12}(x_{2}-x_{1})(v\otimes u\otimes
w)=Y(x_{2})(1\otimes Y(x_{1}))T^{12}(x_{2}-x_{1})(v\otimes u\otimes
w).
\end{eqnarray*}
Given that $K$ is non-degenerate, we know $\pi_{2}$ is injective. It
follows that $R(x)(v\otimes u)=T(x)(v\otimes u)$. This proves the
uniqueness.

For $u\in U$, by definition we have
$$Y({\bf 1},x_{1})Y(u,x_{2})w\sim Y(x_{2})(1\otimes
Y(x_{1}))R^{12}(x_{2}-x_{1})({\bf 1}\otimes u\otimes w).$$ On the
other hand, as $Y({\bf 1},x)=1$ we have
$$Y({\bf 1},x_{1})Y(u,x_{2})w=Y(u,x_{2})Y({\bf 1},x_{1})w
=Y(x_{2})(1\otimes Y(x_{1}))(u\otimes {\bf 1}\otimes w).$$ Thus
$$Y(x_{2})(1\otimes
Y(x_{1}))R^{12}(x_{2}-x_{1})({\bf 1}\otimes u\otimes w)\sim
Y(x_{2})(1\otimes Y(x_{1}))(u\otimes {\bf 1}\otimes w).$$ As we have
seen above, this equivalence relation implies equality relation
$$Y(x_{2})(1\otimes
Y(x_{1}))R^{12}(x_{2}-x_{1})({\bf 1}\otimes u\otimes w)=
Y(x_{2})(1\otimes Y(x_{1}))(u\otimes {\bf 1}\otimes w).$$ Just as
above, with $\pi_{2}$ being injective it follows that $R({\bf
1}\otimes u)=u\otimes {\bf 1}$. Using a parallel argument we get
that $R(x)(v\otimes {\bf 1})={\bf 1}\otimes v$ for $v\in V$.

Now, let $v\in V,\ u,u'\in U$ and let $w\in K$. We have
\begin{eqnarray*}
&&Y(v,z_{1})Y(u,x_{1})Y(u',x_{2})w\\
&\sim & Y(x_{1})(1\otimes Y(z_{1}))(1\otimes 1\otimes
Y(x_{2}))R^{12}(x_{1}-z_{1})(v\otimes u\otimes u'\otimes w)\\
&\sim&Y(x_{1})(1\otimes Y(x_{2}))(1\otimes 1\otimes
Y(z_{1}))R^{23}(x_{2}-z_{1})R^{12}(x_{1}-z_{1})(v\otimes u\otimes
u'\otimes w).
\end{eqnarray*}
Furthermore, by Lemma \ref{ldef-module} there exists a nonnegative
integer $k$ such that
\begin{eqnarray*}
&&((x_{1}-x_{2})^{k}Y(x_{1})(1\otimes Y(x_{2}))(1\otimes
1\otimes Y(z_{1}))\\
&&\hspace{1cm}\cdot R^{23}(x_{2}-z_{1})R^{12}(x_{1}-z_{1})(v\otimes
u\otimes
u'\otimes w))|_{x_{1}=x_{2}+x_{0}}\\
&=&x_{0}^{k}Y(x_{2})(Y(x_{0})\otimes 1)(1\otimes 1\otimes
Y(z_{1}))R^{23}(x_{2}-z_{1})R^{12}(x_{2}+x_{0}-z_{1})(v\otimes
u\otimes u'\otimes w).
\end{eqnarray*}
 Thus
\begin{eqnarray}\label{eneed-one}
&&\left((x_{1}-x_{2})^{k}Y(v,z_{1})Y(u,x_{1})Y(u',x_{2})w\right)|_{x_{1}=x_{2}+x_{0}}
\nonumber\\
&\sim&x_{0}^{k}Y(x_{2})(Y(x_{0})\otimes 1)(1\otimes 1\otimes
Y(z_{1}))R^{23}(x_{2}-z_{1})R^{12}(x_{2}+x_{0}-z_{1})(v\otimes
u\otimes u'\otimes w)\nonumber\\
&=&x_{0}^{k}Y(x_{2})(1\otimes Y(z_{1}))(Y(x_{0})\otimes 1\otimes
1)R^{23}(x_{2}-z_{1})R^{12}(x_{2}-z_{1}+x_{0})(v\otimes u\otimes
u'\otimes w).\nonumber\\
\end{eqnarray}
On the other hand, by Lemma \ref{ldef-module} there exists a
nonnegative integer $k'$ such that
$$(x_{1}-x_{2})^{k'}Y(u,x_{1})Y(u',x_{2})\in \Hom
(K,K((x_{1},x_{2})))$$ and
\begin{eqnarray*}
\left((x_{1}-x_{2})^{k'}Y(u,x_{1})Y(u',x_{2})w\right)|_{x_{1}=x_{2}+x_{0}}
=x_{0}^{k'}Y(Y(u,x_{0})u',x_{2})w.
\end{eqnarray*}
Then
\begin{eqnarray}
&&\left((x_{1}-x_{2})^{k'}Y(v,z_{1})Y(u,x_{1})Y(u',x_{2})w\right)|_{x_{1}=x_{2}+x_{0}}
\nonumber\\
&=& x_{0}^{k'}Y(v,z_{1})Y(Y(u,x_{0})u',x_{2})w\nonumber\\
&=&x_{0}^{k'}Y(z_{1})(1\otimes Y(x_{2}))(1\otimes Y(x_{0})\otimes
1)(v\otimes u\otimes u'\otimes w)\nonumber\\
&\sim&x_{0}^{k'}Y(x_{2})(1\otimes Y(z_{1}))
R^{12}(x_{2}-z_{1})(1\otimes Y(x_{0})\otimes 1)(v\otimes u\otimes
u'\otimes w).
\end{eqnarray}
Combining this with (\ref{eneed-one}) we get
\begin{eqnarray*}
&&Y(x_{2})(1\otimes Y(z_{1})) R^{12}(x_{2}-z_{1})(1\otimes
Y(x_{0})\otimes 1)(v\otimes u\otimes u'\otimes w)\\
&\sim&Y(x_{2})(1\otimes Y(z_{1}))(Y(x_{0})\otimes 1\otimes
1)R^{23}(x_{2}-z_{1})R^{12}(x_{2}-z_{1}+x_{0})(v\otimes u\otimes
u'\otimes w).
\end{eqnarray*}
Because both sides lie in $K((x_{2}))((z_{1}))((x_{0}))$, this
similarity relation implies equality relation. Then, with $\pi_{2}$
injective we have
\begin{eqnarray*}
&&R(x_{2}-z_{1})(1\otimes Y(x_{0}))(v\otimes u\otimes u')\\
&=&(Y(x_{0})\otimes 1)R^{23}(x_{2}-z_{1})
R^{12}(x_{2}-z_{1}+x_{0})(v\otimes u\otimes u').
\end{eqnarray*}
This proves
\begin{eqnarray*}
R(x)(1\otimes Y(x_{0}))= (Y(x_{0})\otimes
1)R^{23}(x)R^{12}(x+x_{0}),
\end{eqnarray*}
confirming (\ref{sy12}). The other condition (\ref{ry12}) can be
proved in the same manner. Therefore, $R$ is a twisting operator for
the ordered pair $(U,V)$.

As for the last assertion, it follows from Theorem \ref{tuniversal}
that $\psi$ is a nonlocal-vertex-algebra homomorphism, and we have
$\psi|_{U}=1$ and $\psi|_{V}=1$. Since $K$ is generated by $U\cup
V$, it follows that $\psi$ is onto. For $u\in U,\ v\in V$, we have
$Y(u,x)v\in K[[x]]$ and $[\D,Y(u,x)]=\frac{d}{dx}Y(u,x)$, which
imply $Y(u,x)v=e^{x\D}u_{-1}v$. From this we get $\ker \psi\subset
\ker Y(x)$. Now we show that $\ker Y(x)=0$. For $a,b\in K$, by weak
associativity, we have
$$(x_{0}+x_{2})^{l}Y(a,x_{0}+x_{2})Y(b,x_{2}){\bf 1}
=(x_{0}+x_{2})^{l}Y(Y(a,x_{0})b,x_{2}){\bf 1}$$ for some nonnegative
integer $l$. As $$Y(a,x_{0}+x_{2})Y(b,x_{2}){\bf 1},\ \
Y(Y(a,x_{0})b,x_{2}){\bf 1}\in K((x_{0}))[[x_{2}]],$$ we must have
$$Y(a,x_{0}+x_{2})Y(b,x_{2}){\bf 1}=Y(Y(a,x_{0})b,x_{2}){\bf 1}.$$
{}From this we have $\ker Y(x)\subset \ker Z_{2}$. With $K$
non-degenerate, we have $\ker Z_{2}=0$, so that $\ker Y(x)=0$ and
$\ker \psi=0$. Thus $\psi$ is also injective. Therefore, $\psi$ is
an isomorphism.
\end{proof}

Note that as $U$ and $V$ are subalgebras of $U\otimes_{R}V$, every
$U\otimes_{R}V$-module is naturally a $U$-module and a $V$-module.
Furthermore, we have:

\bl{lURV-module} Let $(W,Y_{W})$ be a $U\otimes_{R}V$-module and let
 $u\in U,\ v\in V$. Then
\begin{eqnarray}\label{eno-pole}
Y_{W}(u,x_{1})Y_{W}(v,x_{2})\in \Hom (W,W((x_{1},x_{2}))),
 \end{eqnarray}
and there exists a nonnegative integer $k$ such that
\begin{eqnarray}
&&(x_{2}-x_{1})^{k}Y_{W}(v,x_{1})Y_{W}(u,x_{2})w\nonumber\\
&=&(x_{2}-x_{1})^{k}Y_{W}(x_{2})(1\otimes
Y_{W}(x_{1}))R^{12}(x_{2}-x_{1})(v\otimes u\otimes w)
\end{eqnarray}
for every $w\in W$. If $R$ is invertible, we also have
\begin{eqnarray}\label{eYWuv-comm}
Y_{W}(u,x_{1})Y_{W}(v,x_{2})w =Y_{W}(x_{2})(1\otimes
Y_{W}(x_{1}))(R^{-1})^{12}(x_{2}-x_{1})(u\otimes v\otimes w).
\end{eqnarray}
\el

\begin{proof} Let $u\in U,\ v\in V$ and let $w\in W$. There exists
$l\in \N$ such that
$$(x_{0}+x_{2})^{l}Y_{W}(u,x_{0}+x_{2})Y_{W}(v,x_{2})w
=(x_{0}+x_{2})^{l}Y_{W}(Y_{R}(u,x_{0})v,x_{2})w.$$ As
$Y_{R}(u,x_{0})v\in (U\otimes_{R}V)[[x_{0}]]$, the expression on the
right side lies in $W((x_{2}))[[x_{0}]]$. This forces the expression
on the left side to lie in $W((x_{2}))[[x_{0}]]\cap
W((x_{0}))((x_{2}))$. Thus
$$(x_{0}+x_{2})^{l}Y_{W}(u,x_{0}+x_{2})Y_{W}(v,x_{2})w\in
W[[x_{0},x_{2}]][x_{2}^{-1}].$$ Applying $e^{-x_{2}\partial/\partial
x_{0}}$ we get
$$x_{0}^{l}Y_{W}(u,x_{0})Y_{W}(v,x_{2})w\in
W[[x_{0},x_{2}]][x_{2}^{-1}],$$ which implies
$Y_{W}(u,x_{0})Y_{W}(v,x_{2})w\in W((x_{0},x_{2}))$. This proves
(\ref{eno-pole}).

For $u\in U,\ v\in V$, skew symmetry relation (\ref{evu-symmetry})
holds. Then the second assertion follows immediately from
\cite{li-qva1} (Proposition 5.2).

Assume that $R$ is invertible. Then (\ref{eyruv-comm}) (in
Proposition \ref{prelations-prod}) holds. By Corollary 5.4 of
\cite{li-qva1}, there exists a nonnegative integer $k$ such that
\begin{eqnarray*}
&&(x_{2}-x_{1})^{k}Y_{W}(u,x_{1})Y_{W}(v,x_{2})w\nonumber\\
&=&(x_{2}-x_{1})^{k}Y_{W}(x_{2})(1\otimes
Y_{W}(x_{1}))(R^{-1})^{12}(x_{2}-x_{1})(u\otimes v\otimes w),
\end{eqnarray*}
which together with weak associativity gives
\begin{eqnarray*}
&&x_{0}^{-1}\delta\left(\frac{x_{1}-x_{2}}{x_{0}}\right)Y_{W}(u,x_{1})Y_{W}(v,x_{2})w\\
&& \ \ \
-x_{0}^{-1}\delta\left(\frac{x_{2}-x_{1}}{-x_{0}}\right)Y_{W}(x_{2})(1\otimes
Y_{W}(x_{1}))(R^{-1})^{12}(x_{2}-x_{1})(u\otimes v\otimes w)\\
&=&x_{1}^{-1}\delta\left(\frac{x_{2}+x_{0}}{x_{1}}\right)Y_{W}(Y_{R}(u,x_{0})v,x_{2})w.
\end{eqnarray*}
Since $Y_{R}(u,x_{0})v\in (U\otimes_{R}V)[[x_{0}]]$, applying
$\Res_{x_{0}}$ we obtain (\ref{eYWuv-comm}).
\end{proof}

On the other hand, we show that a $U$-module structure together with
a compatible $V$-module structure on the same space gives rise to a
$U\otimes_{R}V$-module structure.

\bp{pmodule-construction} Let $U$ and $V$ be nonlocal vertex
algebras and let $R(x)$ be an invertible twisting operator for
$(U,V)$. Let $W$ be a vector space equipped with a $U$-module
structure $(W,Y_{W}^{U})$ and a $V$-module structure
$(W,Y_{W}^{V})$. Assume that for $u\in U,\ v\in V$,
\begin{eqnarray}\label{ezero-compat}
Y_{W}^{U}(u,x_{1})Y_{W}^{V}(v,x_{2})\in \Hom (W,W((x_{1},x_{2}))),
\end{eqnarray}
(\ref{eYWuv-comm}) holds, and there exists a nonnegative integer $k$
such that
\begin{eqnarray}\label{ekvu=kuv}
&&(x_{2}-x_{1})^{k}Y_{W}^{V}(v,x_{1})Y_{W}^{U}(u,x_{2})w\nonumber\\
&=&(x_{2}-x_{1})^{k}Y_{W}^{U}(x_{2})(1\otimes
Y_{W}^{V}(x_{1}))R^{12}(x_{2}-x_{1})(v\otimes u\otimes w)
\end{eqnarray}
for every $w\in W$. Then there exists a module structure $Y^{W}_{R}$
on $W$ for $U\otimes_{R}V$, extending $Y_{W}^{U}$ and $Y_{W}^{V}$
uniquely. \ep

\begin{proof} For $u\in U,\ v\in V$, we define $Y_{R}^{W}(u\otimes v,x)\in (\End
W)[[x,x^{-1}]]$ by
\begin{eqnarray*}
Y_{R}^{W}(u\otimes v,x)w
=\left(Y_{W}^{U}(u,x_{1})Y_{W}^{V}(v,x)w\right)|_{x_{1}=x}
\end{eqnarray*}
for $w\in W$. Notice that since
$Y_{W}^{U}(u,x_{1})Y_{W}^{V}(v,x)w\in W((x_{1},x))$ by assumption,
$$\left(Y_{W}^{U}(u,x_{1})Y_{W}^{V}(v,x)w\right)|_{x_{1}=x}$$
exists in $W((x))$. Thus $Y_{R}^{W}(u\otimes v, x)$ is well defined
as an element of $\Hom (W,W((x)))$. It is clear that $Y_{R}^{W}({\bf
1}\otimes {\bf 1},x)=1_{W}$.

To establish weak associativity, let $u,u'\in U,\ v,v'\in V$ and let
$w\in W$. We shall use Lemma \ref{ldef-module}. By definition we
have
\begin{eqnarray}\label{efirst-expression}
&&Y_{R}^{W}(u\otimes v,x_{1})Y_{R}^{W}(u'\otimes v',x_{2})w\nonumber\\
&=&\left(Y_{W}^{U}(u,z_{1})Y_{W}^{V}(v,x_{1})Y_{W}^{U}(u',z_{2})
Y_{W}^{V}(v',x_{2})w\right)|_{z_{1}=x_{1},\; z_{2}=x_{2}}.
\end{eqnarray}
We next show that there exists a nonnegative integer $k$ independent
of $w$ such that
\begin{eqnarray}\label{pcondition}
p(z_{1},x_{1},z_{2},x_{2})^{k}Y_{W}^{U}(u,z_{1})Y_{W}^{V}(v,x_{1})Y_{W}^{U}(u',z_{2})
Y_{W}^{V}(v',x_{2})w \in W((z_{1},x_{1},z_{2},x_{2})),
\end{eqnarray}
where
$$p(z_{1},x_{1},z_{2},x_{2})
=(z_{1}-z_{2})(x_{1}-z_{2})(x_{1}-x_{2}).$$
 {}From assumption we have
\begin{eqnarray}\label{every-first}
Y_{W}^{U}(u,z_{1})Y_{W}^{V}(v,x_{1})Y_{W}^{U}(u',z_{2})
Y_{W}^{V}(v',x_{2})w\in W((z_{1},x_{1}))((z_{2},x_{2})).
\end{eqnarray}

Let $k$ be a nonnegative integer independent of $w$ such that
(\ref{ekvu=kuv}) with $u'$ in place of $u$ holds and such that
$x^{k}R(x)(v\otimes u')\in U\otimes V\otimes \C[[x]]$. In view of
Lemma \ref{ldef-module} for $(W,Y_{W}^{V})$, we may also assume that
\begin{eqnarray*}
&&(x_{1}-x_{2})^{k}(1\otimes 1\otimes Y_{W}^{V}(x_{1}))(1\otimes
1\otimes 1\otimes Y_{W}^{V}(x_{2}))R^{23}(\xi)(u\otimes v\otimes
u'\otimes v'\otimes w)\\
&&\ \ \in U\otimes U\otimes W((x_{1},x_{2}))\otimes \C((\xi)).
\end{eqnarray*}
Then
\begin{eqnarray}\label{efirst-1}
&&(x_{1}-x_{2})^{k}(z_{2}-x_{1})^{k}Y_{W}^{U}(u,z_{1})Y_{W}^{V}(v,x_{1})Y_{W}^{U}(u',z_{2})
Y_{W}^{V}(v',x_{2})w\nonumber\\
 &=&(x_{1}-x_{2})^{k}(Y_{W}^{U}(z_{1})(1\otimes Y_{W}^{U}(z_{2}))(1\otimes 1\otimes
Y_{W}^{V}(x_{1}))(1\otimes 1\otimes 1\otimes Y_{W}^{V}(x_{2}))\nonumber\\
&&\ \ \ \ \cdot (z_{2}-x_{1})^{k}R^{23}(z_{2}-x_{1})(u\otimes
v\otimes u'\otimes v'\otimes w)\nonumber\\
&\in&W((z_{1}))((z_{2},x_{1},x_{2}))
\end{eqnarray}
(recall (\ref{every-first})). {}From assumption (\ref{eYWuv-comm})
we also have
\begin{eqnarray}\label{esecond-1}
&&Y_{W}^{U}(u,z_{1})Y_{W}^{V}(v,x_{1})Y_{W}^{U}(u',z_{2})
Y_{W}^{V}(v',x_{2})w\nonumber\\
&=&(Y_{W}^{V}(x_{1})(1\otimes Y_{W}^{U}(z_{1}))(1\otimes 1\otimes
Y_{W}^{U}(z_{2}))(1\otimes 1\otimes 1\otimes Y_{W}^{V}(x_{2}))\nonumber\\
&&(R^{-1})^{12}(x_{1}-z_{1})(u\otimes v\otimes u'\otimes v'\otimes
w)\nonumber\\
&=&(Y_{W}^{V}(x_{1})(1\otimes Y_{W}^{U}(z_{1}))(1\otimes 1\otimes
Y_{W}^{V}(x_{2}))(1\otimes 1\otimes 1\otimes Y_{W}^{U}(z_{2}))\nonumber\\
&&(R^{-1})^{34}(x_{2}-z_{2})(R^{-1})^{12}(x_{1}-z_{1})(u\otimes
v\otimes u'\otimes v'\otimes
w)\nonumber\\
&=&(Y_{W}^{V}(x_{1})(1\otimes Y_{W}^{V}(x_{2}))(1\otimes 1\otimes
Y_{W}^{U}(z_{1}))(1\otimes 1\otimes 1\otimes Y_{W}^{U}(z_{2}))\nonumber\\
&&(R^{-1})^{23}(x_{2}-z_{1})(R^{-1})^{34}(x_{2}-z_{2})(R^{-1})^{12}(x_{1}-z_{1})(u\otimes
v\otimes u'\otimes v'\otimes w).\ \ \ \
\end{eqnarray}
We may choose $k$ (independent of $w$) so large that
\begin{eqnarray*}
&&(z_{1}-z_{2})^{k}(1\otimes 1\otimes Y^{U}_{W}(z_{1}))(1\otimes
1\otimes 1\otimes Y^{U}_{W}(z_{2}))\\
&&\cdot
(R^{-1})^{23}(\xi_{1})(R^{-1})^{34}(\xi_{2})(R^{-1})^{12}(\xi_{3})(u\otimes
v\otimes u'\otimes v')\\
&&\ \ \in V\otimes V\otimes \C((\xi_{1},\xi_{2},\xi_{3}))\otimes
W((z_{1},z_{2})).
\end{eqnarray*}
Then from (\ref{esecond-1}) we see that
\begin{eqnarray*}
(z_{1}-z_{2})^{k}Y_{W}^{U}(u,z_{1})Y_{W}^{V}(v,x_{1})Y_{W}^{U}(u',z_{2})
Y_{W}^{V}(v',x_{2})w
\end{eqnarray*}
lies in $W((x_{1}))((x_{2},z_{1},z_{2}))$. Combining this with
(\ref{efirst-1}) we obtain (\ref{pcondition}).

Using (\ref{efirst-expression}) and (\ref{pcondition}) we get
\begin{eqnarray*}
&&(x_{1}-x_{2})^{3k}Y_{R}^{W}(u\otimes
v,x_{1})Y_{R}^{W}(u'\otimes v',x_{2})w\\
&=&\left(p(z_{1},x_{1},z_{2},x_{2})^{k}Y_{W}^{U}(u,z_{1})Y_{W}^{V}(v,x_{1})Y_{W}^{U}(u',z_{2})
Y_{W}^{V}(v',x_{2})w
\right)|_{z_{1}=x_{1},z_{2}=x_{2}}\\
 && \in W((x_{1},x_{2})),
\end{eqnarray*}
and we have
\begin{eqnarray*}
&&\left(p(z_{1},x_{1},z_{2},x_{2})^{k}Y_{W}^{U}(u,z_{1})Y_{W}^{V}(v,x_{1})Y_{W}^{U}(u',z_{2})
Y_{W}^{V}(v',x_{2})w\right)|_{z_{1}=z_{2}+x_{0},\; x_{1}=x_{2}+x_{0}}\\
 &=&(z_{1}-z_{2})^{k}(x_{1}-x_{2})^{k}(Y_{W}^{U}(z_{1})(1\otimes Y_{W}^{U}(z_{2}))(1\otimes 1\otimes
Y_{W}^{V}(x_{1}))(1\otimes 1\otimes 1\otimes Y_{W}^{V}(x_{2}))\\
&&\ \ \ \ \cdot (z_{2}-x_{1})^{k}R^{23}(z_{2}-x_{1})(u\otimes
v\otimes u'\otimes v'\otimes w))|_{z_{1}=z_{2}+x_{0},\; x_{1}=x_{2}+x_{0}}\\
&=&(z_{1}-z_{2})^{k}(x_{1}-x_{2})^{k}Y_{W}^{U}(z_{1})(1\otimes
Y_{W}^{U}(z_{2}))(1\otimes 1\otimes
Y_{W}^{V}(x_{1}))(1\otimes 1\otimes 1\otimes Y_{W}^{V}(x_{2}))\\
&&\ \ \ \ \cdot (x_{1}-z_{2})^{k}R^{23}(-x_{1}+z_{2})(u\otimes
v\otimes u'\otimes v'\otimes w)|_{z_{1}=z_{2}+x_{0},\; x_{1}=x_{2}+x_{0}}\\
&=&x_{0}^{2k}Y_{W}^{U}(z_{2})(1\otimes
Y_{W}^{V}(x_{2}))(Y(x_{0})\otimes Y(x_{0})\otimes
1)\\
&&\ \ \ \ \cdot
(x_{0}-z_{2}+x_{2})^{k}R^{23}(-x_{0}+z_{2}-x_{2})(u\otimes v\otimes
u'\otimes v'\otimes w).
\end{eqnarray*}
Using substitution $z_{2}=x_{2}$ we get
\begin{eqnarray*}
&&\left((x_{1}-x_{2})^{3k}Y_{R}^{W}(u\otimes
v,x_{1})Y_{R}^{W}(u'\otimes v',x_{2})w\right)|_{x_{1}=x_{2}+x_{0}}\\
&=&x_{0}^{3k}Y_{W}^{U}(z_{2})(1\otimes
Y_{W}^{V}(x_{2}))(Y(x_{0})\otimes Y(x_{0})\otimes
1)\\
&&\ \ \ \ \cdot R^{23}(-x_{0})(u\otimes v\otimes u'\otimes v'\otimes
w))|_{z_{2}=x_{2}}.
\end{eqnarray*}
 On the other hand, we have
\begin{eqnarray*}
&&Y_{R}^{W}\left(Y_{R}(u\otimes v,x_{0})(u'\otimes v'),x_{2}\right)w\\
&=&Y_{R}^{W}(x_{2})(Y(x_{0})\otimes Y(x_{0})\otimes 1)
R^{23}(-x_{0})(u\otimes v\otimes u'\otimes v'\otimes w)\\
&=&Y_{W}^{U}(z_{2})(1\otimes Y_{W}^{V}(x_{2}))(Y(x_{0})\otimes
Y(x_{0})\otimes 1) R^{23}(-x_{0})(u\otimes v\otimes u'\otimes
v'\otimes w)|_{z_{2}=x_{2}}.
\end{eqnarray*}
Combining the last two equations we obtain the desired weak
associativity relation. Therefore, $(W,Y_{W}^{R})$ carries the
structure of a $U\otimes_{R}V$-module. It is clear that the linear
map $Y_{W}^{R}$ extends both $Y_{W}^{U}$ and $Y_{W}^{V}$.

Suppose that $\bar{Y}_{W}$ is another $U\otimes_{R}V$-module
structure on $W$, which also extends both $Y_{W}^{U}$ and
$Y_{W}^{V}$. For $u\in U,\ v\in V$, with the assumption
(\ref{ezero-compat}), by Lemma \ref{ldef-module} we have
\begin{eqnarray}
\bar{Y}_{W}(Y_{R}(u,x_{0})v,x_{2})
=\left(\bar{Y}_{W}(u,x_{1})\bar{Y}_{W}(v,x_{2})\right)|_{x_{1}=x_{2}+x_{0}}.
\end{eqnarray}
Recall that $Y_{R}(u,x)v\in (U\otimes_{R}V)[[x]]$ and $$u\otimes
v=u_{-1}v=\lim_{x\rightarrow 0}Y_{R}(u,x)v.$$ Then
$$\bar{Y}_{W}(u\otimes v,x)=\lim_{x_{0}\rightarrow
0}\bar{Y}_{W}(Y_{R}(u,x_{0})v,x)=\left(\bar{Y}_{W}(u,x_{1})\bar{Y}_{W}(v,x)\right)|_{x_{1}=x}
=Y_{R}^{W}(u\otimes v,x).$$ Thus we have $\bar{Y}_{W}=Y_{R}^{W}$,
proving the uniqueness assertion.
\end{proof}

\section{Twisted tensor product of quantum vertex algebras}
In this section we study twisted tensor product $U\otimes_{R}V$ with
$U$ and $V$ (weak) quantum vertex algebras.

First we recall the notion of weak quantum vertex algebra from
\cite{li-qva1}.

\bd{dwqva} {\em A {\em weak quantum vertex algebra} is a nonlocal
vertex algebra $V$ which satisfies \emph{$\S$-locality} in the sense
that for $u, v\in V$, there exist
\begin{eqnarray*}
u^{(i)}, \ v^{(i)}\in V,\ \ f_{i}(x)\in  \C((x)) \ \ (i=1,\dots,r)
\end{eqnarray*}
(finitely many) such that
\begin{eqnarray}\label{es-locality}
(x_{1}-x_{2})^{k}Y(u,x_{1})Y(v,x_{2})
=(x_{1}-x_{2})^{k}\sum_{i=1}^{r}f_{i}(x_{2}-x_{1})Y(v^{(i)},x_{2})Y(u^{(i)},x_{1})
\end{eqnarray}
for some nonnegative integer $k$.} \ed

The following basic facts can be found in \cite{li-qva1}:

\bp{ps-jacobi} Let $V$ be a nonlocal vertex algebra and
let
\begin{eqnarray*}
u,v,\ u^{(i)}, \ v^{(i)}\in V,\ \ f_{i}(x)\in  \C((x)) \ \
(i=1,\dots,r).
\end{eqnarray*}
Then the $\S$-locality relation (\ref{es-locality}) is equivalent to
\begin{eqnarray}\label{esjacobi}
&&x_{0}^{-1}\delta\left(\frac{x_{1}-x_{2}}{x_{0}}\right)Y(u,x_{1})Y(v,x_{2})
\nonumber\\
&&\hspace{2cm}-x_{0}^{-1}\delta\left(\frac{x_{2}-x_{1}}{-x_{0}}\right)\sum_{i=1}^{r}
f_{i}(-x_{0})Y(v^{(i)},x_{2})Y(u^{(i)},x_{1})\nonumber\\
&&=x_{2}^{-1}\delta\left(\frac{x_{1}-x_{0}}{x_{2}}\right)Y(Y(u,x_{0})v,x_{2})
\end{eqnarray}
(the {\em $\S$-Jacobi identity}), and it is also equivalent to
\begin{eqnarray}\label{es-skewsymmetry}
Y(u,x)v=e^{x\D}\sum_{i=1}^{r}f_{i}(-x)Y(v^{(i)},-x)u^{(i)}
\end{eqnarray}
(the {\em $\S$-skew symmetry}). \ep

\bp{pmodule-wqva} Let $V$ be a weak quantum vertex algebra and let
$(W,Y_{W})$ be a module for $V$ viewed as a nonlocal vertex algebra.
Assume
\begin{eqnarray*}
u,v,\ u^{(i)}, \ v^{(i)}\in V,\ \ f_{i}(x)\in  \C((x)) \ \
(i=1,\dots,r)
\end{eqnarray*}
such that (\ref{es-skewsymmetry}) holds. Then
\begin{eqnarray*}
&&x_{0}^{-1}\delta\left(\frac{x_{1}-x_{2}}{x_{0}}\right)Y_{W}(u,x_{1})Y_{W}(v,x_{2})
\nonumber\\
&&\hspace{2cm}-x_{0}^{-1}\delta\left(\frac{x_{2}-x_{1}}{-x_{0}}\right)\sum_{i=1}^{r}
f_{i}(-x_{0})Y_{W}(v^{(i)},x_{2})Y_{W}(u^{(i)},x_{1})\nonumber\\
&&=x_{2}^{-1}\delta\left(\frac{x_{1}-x_{0}}{x_{2}}\right)Y_{W}(Y(u,x_{0})v,x_{2}).
\end{eqnarray*}
\ep

\br{rweak-smap} {\em  {}From definition, a nonlocal vertex algebra
$V$ is a weak quantum vertex algebra if and only if there exists a
linear map
$$\S(x): V\otimes V\rightarrow V\otimes V\otimes \C((x))$$
satisfying the condition that for $u,v\in V$, there exists $k\in \N$
such that
\begin{eqnarray}
&&(x_{1}-x_{2})^{k}Y(x_{1})(1\otimes Y(x_{2}))(u\otimes v\otimes
w)\nonumber\\
&=&(x_{1}-x_{2})^{k}Y(x_{2})(1\otimes
Y(x_{1}))\S^{12}(x_{2}-x_{1})(v\otimes u\otimes w)
\end{eqnarray}
for every $w\in V$, or equivalently,
\begin{eqnarray}
Y(x)(u\otimes v)=e^{x\D}Y(-x)\S(-x)(v\otimes u).
\end{eqnarray}}
\er

A {\em rational quantum Yang-Baxter operator} on a vector space $U$
is a linear operator
$$\S(x):\ U\otimes U\rightarrow U\otimes U\otimes \C((x))$$
satisfying the quantum Yang-Baxter equation
$$\S^{12}(x)\S^{13}(x+z)\S^{23}(z)=\S^{23}(z)\S^{13}(x+z)\S^{12}(x).$$
It is said to be {\em unitary} if
$$\S(x)\S^{21}(-x)=1,$$
where $\S^{21}(x)=\sigma \S(x)\sigma$ with $\sigma$ denoting the
flip operator on $U\otimes U$.

\bd{dqva} {\em A {\em quantum vertex algebra} is a weak quantum
vertex algebra $V$ equipped with a unitary rational quantum
Yang-Baxter operator $\S(x)$ on $V$, satisfying
\begin{eqnarray}
&&\S(x)({\bf 1}\otimes v)={\bf 1}\otimes v\ \ \ \mbox{ for }v\in
V,\label{esvacuum}\\
&&[\D\otimes 1, \S(x)]=-\frac{d}{dx}\S(x),\label{d1s}\\
&&Y(u,x)v=e^{x\D}Y(-x)\S(-x)(v\otimes u)\ \ \mbox{ for }u,v\in V,\\
&&\S(x_{1})(Y(x_{2})\otimes 1)=(Y(x_{2})\otimes
1)\S^{23}(x_{1})\S^{13}(x_{1}+x_{2}).\label{sy1}
\end{eqnarray}
We denote a quantum vertex algebra by a pair $(V,\S)$.} \ed

Note that this very notion is a slight modification of the same
named notion in \cite{li-qva1} and \cite{li-qva2} with extra axioms
(\ref{esvacuum}) and (\ref{sy1}).

The following are some axiomatic results:

\bl{vacuum-ek} Let $V$ be a nonlocal vertex algebra and let $\S(x)$
be a unitary rational quantum Yang-Baxter equation on $V$. Then
(\ref{esvacuum}), (\ref{d1s}), and (\ref{sy1}) are equivalent to
\begin{eqnarray}
&&\S(x)(v\otimes {\bf 1})=v\otimes {\bf 1}\ \ \ \mbox{ for }v\in
V,\label{esv1=v1}\\
&&[1\otimes \D,\S^{-1}(x)]=\frac{d}{dx}\S^{-1}(x),\\
&&\S(x_{1})(1\otimes Y(x_{2}))=(1\otimes
Y(x_{2}))\S^{12}(x_{1}-x_{2})\S^{13}(x_{1}),\label{psy2}
\end{eqnarray}
respectively.\el

\begin{proof} By unitarity
we have $\S^{-1}(x)=\S^{21}(-x)=\sigma \S(-x)\sigma$.  For $v\in V$,
we have
\begin{eqnarray*}
&&\S(x)(v\otimes {\bf 1})=\sigma\sigma \S(x)\sigma({\bf 1}\otimes
v)=\sigma \S^{-1}(-x)({\bf 1}\otimes v),\\
&&\S(x)({\bf 1}\otimes v)=\sigma\sigma \S(x)\sigma(v\otimes {\bf
1})=\sigma \S^{-1}(-x)(v\otimes {\bf 1}).
\end{eqnarray*}
It follows that $\S(x)({\bf 1}\otimes v)={\bf 1}\otimes v$ if and
only if $\S(x)(v\otimes {\bf 1})=v\otimes {\bf 1}$.

Assuming (\ref{d1s}), that is $[\D\otimes
1,\S(x)]=-\frac{d}{dx}\S(x)$, we have
\begin{eqnarray*}
[1\otimes \D, \S^{-1}(x)]=\sigma [\D\otimes 1,\sigma
\S^{-1}(x)\sigma]\sigma=\sigma [\D\otimes 1, \S(-x)]\sigma
=\frac{d}{dx}\S^{-1}(x).
\end{eqnarray*}
Similarly, assuming $[1\otimes \D,
\S^{-1}(x)]=\frac{d}{dx}\S^{-1}(x)$ we have
\begin{eqnarray*}
[\D\otimes 1, \S(x)]=\sigma [1\otimes \D,\sigma
\S(x)\sigma]\sigma=\sigma [1\otimes \D, \S^{-1}(-x)]\sigma
=-\frac{d}{dx}\S(x).
\end{eqnarray*}

Note that (\ref{sy1}) amounts to
$$\S^{-1}(x_{1})(Y(x_{2})\otimes 1)
=(Y(x_{2})\otimes
1)(\S^{-1})^{13}(x_{1}+x_{2})(\S^{-1})^{23}(x_{1}),$$ which is
$$\sigma \S(-x_{1})\sigma(Y(x_{2})\otimes 1)
=(Y(x_{2})\otimes
1)\sigma^{13}\S^{13}(-x_{1}-x_{2})\sigma^{13}\sigma^{23}\S^{23}(-x_{1})\sigma^{23}.$$
The latter amounts to
\begin{eqnarray}\label{emiddle}
\S(-x_{1})\sigma(Y(x_{2})\otimes 1) =\sigma(Y(x_{2})\otimes
1)\sigma^{13}\S^{13}(-x_{1}-x_{2})\sigma^{13}\sigma^{23}\S^{23}(-x_{1})\sigma^{23}.
\end{eqnarray}
As
$$\sigma (Y(x)\otimes 1)=(1\otimes
Y(x))\sigma^{12}\sigma^{23},$$ (\ref{emiddle}) amounts to
\begin{eqnarray*}
&&\S(-x_{1})(1\otimes Y(x_{2}))\\
&=&(1\otimes
Y(x_{2}))\sigma^{12}\sigma^{23}\sigma^{13}\S^{13}(-x_{1}-x_{2})
\sigma^{13}\sigma^{23}\S^{23}(-x_{1})\sigma^{12}\\
&=&(1\otimes Y(x_{2}))\sigma^{23}\S^{13}(-x_{1}-x_{2})
\sigma^{23}\sigma^{12}\S^{23}(-x_{1})\sigma^{12}\\
&=&(1\otimes Y(x_{2}))\S^{12}(-x_{1}-x_{2})\S^{13}(-x_{1}),
\end{eqnarray*}
which is a version of (\ref{psy2}).
\end{proof}

\bl{lqva-twisting} Let $(V,\S)$ be a quantum vertex algebra. Set
$$R(x)=\S(x)\sigma:\ V\otimes V\rightarrow V\otimes V\otimes \C((x)).$$
Then $R(x)$ is an invertible twisting operator for the ordered pair
$(V,V)$. \el

\begin{proof} As $\S(x)$ is unitary, it is clear that $R(x)$ is invertible.
{}From Lemma \ref{vacuum-ek} we have
\begin{eqnarray*}
&&R(x)({\bf 1}\otimes u)=\S(x)(u\otimes {\bf 1})=u\otimes {\bf 1},\\
&&R(x)(v\otimes {\bf 1})=\S(x)({\bf 1}\otimes v)={\bf 1}\otimes v
\end{eqnarray*}
for $u,v\in V$. Notice that
\begin{eqnarray*}
&&\sigma (Y(x)\otimes 1)=(1\otimes Y(x))\sigma^{12}\sigma^{23},\\
&&\sigma (1\otimes Y(x))=(Y(x)\otimes 1)\sigma^{23}\sigma^{12}.
\end{eqnarray*}
Using this, (\ref{sy1}), and (\ref{psy2}), we obtain (\ref{ry12})
and (\ref{sy12}).
\end{proof}

The following is straightforward (cf. \cite{ek}, \cite{li-qva1}):

\bp{pnon-degenerate} Let $V$ be a weak quantum vertex algebra.
Assume that $V$ is non-degenerate. Then there exists a linear map
$\S(x): V\otimes V\rightarrow V\otimes V\otimes \C((x))$, which is
uniquely determined by
\begin{eqnarray*}
Y(u,x)v=e^{x\D}Y(-x)\S(-x)(v\otimes u) \ \ \ \mbox{for }u,v\in V.
\end{eqnarray*}
Furthermore, $(V,\S)$ carries the structure of a quantum vertex
algebra and the following relation holds
\begin{eqnarray}
[1\otimes \D, \S(x)]=\frac{d}{dx}\S(x).
\end{eqnarray}
 \ep

In view of Proposition \ref{pnon-degenerate}, the term
``non-degenerate quantum vertex algebra'' without referring a
quantum Yang-Baxter operator is non-ambiguous.

 \bp{pweak-qva} Let $U$ and $V$ be weak quantum vertex
algebras and let $R(x)$ be an invertible twisting operator for
$(U,V)$. Then $U\otimes_{R}V$ is a weak quantum vertex algebra. \ep

\begin{proof} From Remark \ref{rweak-smap}, there are linear maps
$$\S_{U}(x): U\otimes U\rightarrow U\otimes U\otimes \C((x))\ \
\mbox{and }\  \S_{V}(x): V\otimes V\rightarrow V\otimes V\otimes
\C((x))$$ such that for $u, u'\in U,\ v, v'\in V$,
\begin{eqnarray*}
&&Y_{U}(x)(u\otimes u')=e^{x\D}Y_{U}(-x)\S_{U}(-x)(u'\otimes
u),\\
&&Y_{V}(x)(v\otimes v')=e^{x\D}Y_{V}(-x)\S_{V}(-x)(v'\otimes v).
\end{eqnarray*}
Let $u, u'\in U,\ v, v'\in V$. Using (\ref{dvv}) we get
\begin{eqnarray*}
&&Y_{R}(u\otimes v,x)(u'\otimes v')\nonumber\\
&=&(Y_{U}(x)\otimes Y_{V}(x))R^{23}(-x)(u\otimes v\otimes
u'\otimes v')\nonumber\\
&=&(e^{x\D}Y_{U}(-x)\otimes
e^{x\D}Y_{V}(-x))\S_{U}^{12}(-x)\sigma^{12}\S_{V}^{34}(-x)\sigma^{34}R^{23}(-x)(u\otimes
v\otimes
u'\otimes v')\nonumber\\
&=&e^{x(\D\otimes 1+1\otimes\D)}(Y_{U}(-x)\otimes
Y_{V}(-x))\S_{U}^{12}(-x)\sigma^{12}\S_{V}^{34}(-x)\sigma^{34}R^{23}(-x)(u\otimes
v\otimes u'\otimes v')\nonumber\\
&=&e^{x(\D\otimes
1+1\otimes\D)}Y_{R}(-x)(R^{-1})^{23}(-x)\S_{U}^{12}(-x)\sigma^{12}\S_{V}^{34}(-x)\sigma^{34}
R^{23}(-x)(u\otimes v\otimes u'\otimes v').
\end{eqnarray*}
Setting
\begin{eqnarray}\label{enew-qybo}
\S_{R}(x)=(R^{-1})^{23}(x)\S_{U}^{12}(x)\sigma^{12}\S_{V}^{34}(x)\sigma^{34}
R^{23}(x)\sigma^{13}\sigma^{24},
\end{eqnarray}
a linear map from $(U\otimes_{R}V)\otimes (U\otimes_{R}V)$ to
$(U\otimes_{R}V)\otimes (U\otimes_{R}V)\otimes \C((x))$, we have
\begin{eqnarray}\label{eskew-twistedproduct}
Y_{R}(u\otimes v,x)(u'\otimes v')=e^{x(\D\otimes
1+1\otimes\D)}Y_{R}(-x)\S_{R}(-x)(u'\otimes v'\otimes u\otimes v).
\end{eqnarray}
By Remark \ref{rweak-smap} or by Proposition \ref{ps-jacobi},
$U\otimes_{R}V$ is a weak quantum vertex algebra.
\end{proof}

Furthermore, we have:

\bp{pqva}  Let $U$ and $V$ be weak quantum vertex algebras and let
$R$ be an invertible twisting operator. Assume that $U$ as a
$U$-module and $V$ as a $V$-module are irreducible and of countable
dimension over $\C$. Then $U\otimes_{R}V$ as a
$U\otimes_{R}V$-module is irreducible. Furthermore, $U\otimes_{R}V$
is a non-degenerate quantum vertex algebra. \ep

\begin{proof} Recall from (\ref{eu-action}) that
$$Y_{R}(u,x)(u'\otimes v)=Y(u,x)u'\otimes v\ \ \mbox{ for }u,u'\in
U,\ v\in V.$$ It follows that for any subspace $A$ of $V$, $U\otimes
A$ is a $U$-submodule of $U\otimes_{R}V$.  As $U$ is an irreducible
$U$-module of countable dimension over $\C$, we have $\Hom
_{U}(U,U)=\C$. Let $P$ be any $U\otimes_{R}V$-submodule of
$U\otimes_{R}V$. Using the $U$-module structure we get $P=U\otimes
A$ for some vector space $A$ of $V$. Since
$${\bf 1}\otimes Y(v,x)a=Y_{R}(v,x)({\bf 1}\otimes a)\in (U\otimes A)[[x,x^{-1}]]$$
for $v\in V,\ a\in A\subset V$, we see that $A$ is a $V$-submodule.
This proves that each $U\otimes_{R}V$-submodule of $U\otimes_{R}V$
is of the form $U\otimes A$ where $A$ is a $V$-submodule of $V$.
Since $V$ is an irreducible $V$-module, $U\otimes_{R}V$ as a
$U\otimes_{R}V$-module is irreducible. It follows from Theorem 3.9
of \cite{li-qva2} that $U\otimes_{R}V$ is non-degenerate. On the
other hand, by Proposition \ref{pweak-qva}, $U\otimes_{R}V$ is a
weak quantum vertex algebra. Therefore, $U\otimes_{R}V$ is a
non-degenerate quantum vertex algebra.
\end{proof}

\section{Smash product of nonlocal vertex algebras}
In this section, we first slightly generalize the smash product
construction of nonlocal vertex algebras, established in
\cite{li-smash}, and we then prove that every  smash product is a
twisted tensor product with respect to a canonical twisting
operator.

We begin by recalling from \cite{li-smash} the basic notions and the
smash product construction. For convenience, we first recall some
necessary classical notions. A {\em coalgebra} is a vector space $C$
(over $\C$) equipped with linear maps
$$\Delta:\ C\rightarrow
C\otimes C \ \mbox{ and } \ \varepsilon:\ C\rightarrow \C,$$
satisfying
\begin{eqnarray*}
&&(1\otimes \Delta)\Delta(b)=(\Delta\otimes 1)\Delta(b),\\
&&(\varepsilon \otimes 1)\Delta(b)=1\otimes b,\ \ (1\otimes
\varepsilon)\Delta(b)=b\otimes 1
\end{eqnarray*}
for $b\in C$. For a coalgebra $C$, a {\em $C$-comodule} is a vector
space $V$ equipped with a linear map $\rho: V\rightarrow C\otimes V$
such that
\begin{eqnarray}
&&(1\otimes \rho)\rho =(\Delta\otimes 1)\rho,\label{ecomodule2}\\
&&(\varepsilon \otimes 1)\rho(v)={\bf 1}\otimes v \ \ \mbox{ for
}v\in V.
\end{eqnarray}

A {\em nonlocal vertex bialgebra} is a nonlocal vertex algebra $H$
equipped with a classical coalgebra structure $(\Delta,\varepsilon)$
such that both $\Delta$ and $\varepsilon$ are homomorphisms of
nonlocal vertex algebras. That is,
\begin{eqnarray}
&&\varepsilon({\bf 1})=1, \ \ \varepsilon
(Y(h,x)h')=\varepsilon(h)\varepsilon(h'),\\
&&\Delta({\bf 1})={\bf 1}\otimes {\bf 1}, \ \ \Delta
(Y(h,x)h')=Y(\Delta(h),x)\Delta(h')\ \ \ \mbox{ for }h,h'\in H.
\end{eqnarray}
A {\em nonlocal vertex $H$-module-algebra} is a nonlocal vertex
algebra $V$ equipped with a module structure for $H$ viewed as a
nonlocal vertex algebra such that
\begin{eqnarray}
&&Y(h,x)v\in V\otimes \C((x))\ (\subset V((x))),\\
&&Y(h,x){\bf 1}=\varepsilon(h){\bf 1},\\
&&Y(h,x)Y(u,z)v=Y(Y(h^{1},x-z)u,z)Y(h^{2},x)v\label{emodulealgebra-3}
\end{eqnarray}
for $h\in H,\ u,v\in V$, where $\Delta(h)=h^{1}\otimes h^{2}$ in the
Sweedler notation. Notice that if $V$ is infinite-dimensional, which
is true most of the time, $V\otimes \C((x))\ne V((x))$.

The following is a simple fact that we shall need later:

\bl{lneed} Let $H$ be a nonlocal vertex bialgebra and let $V$ be a
nonlocal vertex $H$-module-algebra. Then
\begin{eqnarray}
Y(h,z+x)Y(h',z)v=Y(Y(h,x)h',z)v
\end{eqnarray}
for $h,h'\in H,\ v\in V$. \el

\begin{proof} Let $h,h'\in H,\ v\in V$. There exists a nonnegative integer $l$ such
that
\begin{eqnarray*}
(x+z)^{l}Y(h,x+z)Y(h',z)v=(x+z)^{l}Y(Y(h,x)h',z)v.
\end{eqnarray*}
(Note that $Y(h,x+z)Y(h',z)v$ exists in $V((x))((z))$.) As
$Y(h',z)v\in V\otimes \C((z))$, we have
$$Y(h,x_{1})Y(h',z)v\in V((x_{1}))\otimes \C((z)).$$
Then $Y(h,z+x)Y(h',z)v$ exists in $V((z))[[x]]$.
 Replace $l$ with a large one if necessary, so that
$$x_{1}^{l}Y(h,x_{1})Y(h',z)v\in V[[x_{1}]]\otimes \C((z)).$$
Because of this we have
$$(x+z)^{l}Y(h,x+z)Y(h',z)v=(z+x)^{l}Y(h,z+x)Y(h',z)v.$$
Then
\begin{eqnarray*}
(z+x)^{l}Y(h,z+x)Y(h',z)v=(x+z)^{l}Y(Y(h,x)h',z)v.
\end{eqnarray*}
Multiplying both sides by $(z+x)^{-l}\ (\in \C((z))[[x]])$, we
obtain the desired relation.
\end{proof}

\bd{dcomodule-algebra} {\em Let $H$ be a nonlocal vertex bialgebra.
A {\em nonlocal vertex $H$-comodule-algebra} is a nonlocal vertex
algebra $V$ equipped with a (left) comodule structure
$$\rho: V\rightarrow H\otimes V$$
for $H$ viewed as a coalgebra such that $\rho$ is a homomorphism of
nonlocal vertex algebras, that is,
\begin{eqnarray}
&&\rho({\bf 1}) ={\bf 1}\otimes {\bf 1},\label{ecomodule-vacuum}\\
&&\rho(Y(v,x)v')=(Y(x)\otimes Y(x))\sigma^{23}(\rho(v)\otimes
\rho(v'))\ \ \ \mbox{ for }v,v'\in V.\label{ecomodule1}
\end{eqnarray}
} \ed

Notice that $H$ itself is a nonlocal vertex $H$-comodule-algebra
with $\rho=\Delta$. The following is a slight generalization of the
smash product construction of \cite{li-smash} (cf. \cite{kl}):

\bp{psharp-product} Let $H$ be a nonlocal vertex bialgebra, let $U$
be a nonlocal vertex (left) $H$-module-algebra, and let $V$ be a
nonlocal vertex (left) $H$-comodule-algebra. For $u,u'\in U,\
v,v'\in V$, define
\begin{eqnarray}
Y_{\sharp}(u\otimes v,x)(u'\otimes v')=Y(u,x)Y(b^{1}(v),x)u'\otimes
Y(v^{2},x)v',
\end{eqnarray}
where $\rho(v)= b^{1}(v)\otimes v^{2}$. Then $(U\otimes
V,Y_{\sharp},{\bf 1}\otimes {\bf 1})$ carries the structure of a
nonlocal vertex algebra, which we denote by $U\sharp V$. \ep

\begin{proof} For $u,u'\in U,\ v,v'\in V$,
as $Y(b^{1}(v),x)u'\in U\otimes \C((x))$ from definition, we see
that $Y(u,x)Y(b^{1}(v),x)u'\in U((x))$, so that $Y_{\sharp}(u\otimes
v,x)(u'\otimes v')$ is a well defined element of $(U\otimes
V)((x))$. Setting $u={\bf 1},\ v={\bf 1}$, we have
$$Y_{\sharp}({\bf 1}\otimes {\bf 1},x)(u'\otimes v')=Y({\bf 1},x)Y({\bf 1},x)u'\otimes
Y({\bf 1},x)v'=u'\otimes v',$$ as $\rho({\bf 1})={\bf 1}\otimes {\bf
1}$. On the other hand, setting $u'={\bf 1},\ v'={\bf 1}$, we get
\begin{eqnarray*}
Y_{\sharp}(u\otimes v,x)({\bf 1}\otimes {\bf 1})
&=&Y(u,x)Y(b^{1}(v),x){\bf 1}\otimes Y(v^{2},x){\bf 1}\\
&=&\varepsilon(b^{1}(v))Y(u,x){\bf 1}\otimes Y(v^{2},x){\bf 1}\\
&=&Y(u,x){\bf 1}\otimes Y(v,x){\bf 1},
\end{eqnarray*}
which implies
$$Y_{\sharp}(u\otimes v,x)({\bf 1}\otimes {\bf
1})\in (U\otimes V)[[x]]\ \mbox{ and }\ \lim_{x\rightarrow
0}Y_{\sharp}(u\otimes v,x)({\bf 1}\otimes {\bf 1})=u\otimes v.$$ For
weak associativity, let $u,u',u''\in U,\ v,v',v''\in V$. Writing
$$\rho(v)=b^{1}(v)\otimes v^{2},\ \ \ \rho(v')=b^{1}(v')\otimes
v'^{2},$$ we have
\begin{eqnarray*}
&&Y_{\sharp}(u\otimes v,x_{1})Y_{\sharp}(u'\otimes
v',x_{2})(u''\otimes v'')\\
&=&Y(u,x_{1})Y(b^{1}(v),x_{1})Y(u',x_{2})Y(b^{1}(v'),x_{2})u''
\otimes Y(v^{2},x_{1})Y(v'^{2},x_{2})v''\\
&=&Y(u,x_{1})Y(Y(b^{1}(v)^{1},x_{1}-x_{2})u',x_{2})Y(b^{1}(v)^{2},x_{1})Y(b^{1}(v'),x_{2})u''\\
&&\ \ \ \ \ \otimes Y(v^{2},x_{1})Y(v'^{2},x_{2})v'',
\end{eqnarray*}
using (\ref{emodulealgebra-3}), and we have
\begin{eqnarray*}
&&Y_{\sharp}\left(Y_{\sharp}(u\otimes v,x_{0})(u'\otimes
v'),x_{2}\right)(u''\otimes v'')\\
&=&Y_{\sharp}\left(Y(u,x_{0})Y(b^{1}(v),x_{0})u'\otimes
Y(v^{2},x_{0})v',x_{2}\right)(u''\otimes v'')\\
&=&
Y\left(Y(u,x_{0})Y(b^{1}(v),x_{0})u',x_{2}\right)Y(b^{1}(Y(v^{2},x_{0})v'),x_{2})u''\otimes
Y((Y(v^{2},x_{0})v')^{2},x_{2})v''\\
&=&
Y\left(Y(u,x_{0})Y(b^{1}(v),x_{0})u',x_{2}\right)Y(Y(b^{1}(v^{2}),x_{0})b^{1}(v'),x_{2})u''\otimes
Y(Y(v^{22},x_{0})v'^{2}),x_{2})v'',
\end{eqnarray*}
using (\ref{ecomodule1}). Note that (\ref{ecomodule2}) implies
$$b^{1}(v)^{1}\otimes b^{1}(v)^{2}\otimes
v^{2}\otimes b^{1}(v')\otimes v'^{2} =b^{1}(v)\otimes
b^{1}(v^{2})\otimes v^{22}\otimes b^{1}(v')\otimes v'^{2},$$ which
is equivalent to
$$b^{1}(v)^{1}\otimes b^{1}(v)^{2}\otimes
b^{1}(v')\otimes v^{2}\otimes v'^{2} =b^{1}(v)\otimes
b^{1}(v^{2})\otimes b^{1}(v')\otimes v^{22}\otimes v'^{2}.$$ Then,
using Lemma \ref{ldef-module} we obtain the desired weak
associativity relation.
\end{proof}

As the main result of this section we have:

\bp{pequivalence} Let $H,U,V$ be given as in Proposition
\ref{psharp-product}. Define a linear map
$$R(x): V\otimes U\rightarrow U\otimes V\otimes \C((x))$$ by
$$R(x)(v\otimes u)=Y(b^{1}(v),-x)u\otimes v^{2}
\ \ \ \mbox{ for }v\in V,\ u\in U.$$ Then $R(x)$ is a twisting
operator for the pair $(U,V)$. Furthermore, we have $U\sharp V=
U\otimes_{R}V$. \ep

\begin{proof} Let $u\in U,\ v\in V$. As
$\rho({\bf 1})={\bf 1}\otimes {\bf 1},\ \ (\varepsilon\otimes 1)
\rho(v)={\bf 1}\otimes v$, we have
\begin{eqnarray*}
&&R(x)({\bf 1}\otimes u)=Y({\bf 1},-x)u\otimes {\bf 1}=u\otimes {\bf
1},\\
&&R(x)(v\otimes {\bf 1})=Y(b^{1}(v),-x){\bf 1}\otimes v^{2}
=\varepsilon(b^{1}(v)){\bf 1}\otimes v^{2}={\bf 1}\otimes v.
\end{eqnarray*}
On one hand, using (\ref{ecomodule1}) we have
\begin{eqnarray*}
R(z)(Y(x)\otimes 1)(v\otimes v'\otimes u)&=&R(z)(Y(v,x)v'\otimes
u)\\
&=&Y(b^{1}(Y(v,x)v'),-z)u\otimes (Y(v,x)v')^{2}\\
&=&Y\left(Y(b^{1}(v),x)b^{1}(v'),-z\right)u\otimes Y(v^{2},x)v'^{2},
\end{eqnarray*}
noticing that (\ref{ecomodule1}) gives
\begin{eqnarray}
 b^{1}(Y(v,x)v')\otimes
(Y(v,x)v')^{2}=Y(b^{1}(v),x)b^{1}(v')\otimes Y(v^{2},x)v'^{2}.
\end{eqnarray}
On the other hand, using Lemma \ref{lneed} we get
\begin{eqnarray*}
&&(1\otimes Y(x))R^{12}(z-x)R^{23}(z)(v\otimes v'\otimes u)\\
&=&(1\otimes Y(x))R^{12}(z-x)(v\otimes
Y(b^{1}(v'),-z)u\otimes v'^{2})\\
&=&(1\otimes Y(x))\left(Y(b^{1}(v),-z+x)Y(b^{1}(v'),-z)u\otimes
v^{2}\otimes v'^{2}\right)\\
&=&Y(b^{1}(v),-z+x)Y(b^{1}(v'),-z)u\otimes Y(v^{2},x)v'^{2}\\
&=&Y\left(Y(b^{1}(v),x)b^{1}(v'),-z\right)u\otimes Y(v^{2},x)v'^{2}.
\end{eqnarray*}
Consequently, we have
$$R(z)(Y(x)\otimes 1)(v\otimes v'\otimes u)=(1\otimes
Y(x))R^{12}(z-x)R^{23}(z)(v\otimes v'\otimes u).$$ Similarly, we
have
\begin{eqnarray*}
&&R(z)(1\otimes Y(x))(v\otimes u\otimes u')\\
&=&R(z)(v\otimes Y(u,x)u')\\
&=&Y(b^{1}(v),-z)Y(u,x)u'\otimes v^{2}\\
&=&Y(Y(b^{1}(v)^{1},-z-x)u,x)Y(b^{1}(v)^{2},-z)u'\otimes v^{2},
\end{eqnarray*}
 while
\begin{eqnarray*}
&&(Y(x)\otimes 1)R^{23}(z)R^{12}(z+x)(v\otimes u\otimes u')\\
&=&(Y(x)\otimes 1)R^{23}(z)(Y(b^{1}(v),-z-x)u\otimes v^{2}\otimes u')\\
&=&(Y(x)\otimes 1)\left(Y(b^{1}(v),-z-x)u\otimes
Y(b^{1}(v^{2}),-z)u'\otimes
v^{22}\right)\\
&=&Y\left(Y(b^{1}(v),-z-x)u,x\right)Y(b^{1}(v^{2}),-z)u'\otimes
v^{22}.
\end{eqnarray*}
As
$$b^{1}(v)^{1}\otimes b^{1}(v)^{2}\otimes v^{2}=b^{1}(v)\otimes b^{1}(v^{2})\otimes
v^{22}$$ by (\ref{ecomodule2}), we have
$$R(z)(1\otimes Y(x))(v\otimes u\otimes u')=(Y(x)\otimes 1)
R^{23}(z)R^{12}(z+x)(v\otimes u\otimes u').$$ This proves that
$R(x)$ is a twisting operator for the ordered pair $(U,V)$. It
follows immediately from the definitions that $U\sharp V=
U\otimes_{R}V$.
\end{proof}

We note that the smash product nonlocal vertex algebra $U\sharp H$,
which was studied in \cite{li-smash}, is the twisted product
$U\sharp V$ with $V$ specialized to $H$. In \cite{li-smash}, some
interesting examples, related to the vertex operator algebras
associated to even lattices and those associated to
infinite-dimensional Heisenberg algebras, were given.

\end{document}